\documentclass[final,3p,times]{elsarticle}
\usepackage{amssymb}
\usepackage{amsthm}
\usepackage{lineno}
\usepackage{wrapfig}
\usepackage{soul}
\usepackage{footmisc}

\usepackage{multirow}
\usepackage{amsfonts}
\usepackage{amsmath}
\usepackage{url}
\usepackage{todonotes}
\usepackage{units}
\usepackage{physics}
\usepackage{algpseudocode}  % za psevdokodo
\usepackage{algorithm}      % za
\usepackage{cancel}
\usepackage{caption}
\usepackage{booktabs}
\usepackage{varwidth}
\usepackage{subcaption}
\usepackage{xcolor,colortbl}
\definecolor{Gray}{gray}{0.9}
\algnewcommand\algorithmicto{\textbf{to}}
\algnewcommand\algorithmicin{\textbf{in}}
\algnewcommand\algorithmicforeach{\textbf{for each}}
\algrenewtext{For}[3]{\algorithmicfor\ #1 $\gets$ #2\ \algorithmicto\ #3\ \algorithmicdo}
\algdef{S}[FOR]{ForEach}[2]{\algorithmicforeach\ #1\ \algorithmicin\ #2\ \algorithmicdo}

\newcommand{\R}{\mathbb{R}}
\newcommand{\T}{\mathsf{T}}
\renewcommand{\L}{\mathcal{L}}
\renewcommand{\b}{\boldsymbol}

\newcommand{\x}{\b{x}}
\newcommand{\w}{\b{w}}

\newcommand{\lap}{\nabla^2}

\newcommand{\h}{\emph{h}}

\newcommand{\env}{_\mathrm{env}}
\newcommand{\ienv}{_{i,\mathrm{env}}}
\newcommand{\he}{h_\mathrm{e}}
\newcommand*\rot{\rotatebox{90}}

\journal{Journal of Computational Physics}

\begin{document}

\begin{frontmatter}


%    \title{Sharp-interface dendritic envelope growth simulation}
    \title{Meshless interface tracking for the simulation of dendrite envelope growth}

    \author[ijs]{Mitja Jančič\corref{cor1}}
    \ead{mitja.jancic@ijs.si}

    \author[lor,damas]{Miha Založnik}
    \ead{miha.zaloznik@univ-lorraine.fr}
    \author[ijs]{Gregor Kosec}
    \ead{gregor.kosec@ijs.si}

    \address[ijs]{Institute Jožef Stefan, Parallel and Distributed Systems Laboratory, Jamova cesta 39, 1000 Ljubljana, Slovenia}
    % \address[mps]{International Postgraduate School Jožef Stefan, Jamova Cesta 39, 1000 Ljubljana, Slovenia}
    \address[lor]{Universit\'{e} de Lorraine, CNRS, IJL, Nancy, F-54000, France}
    \address[damas]{Laboratory of Excellence DAMAS, Universit\'{e} de Lorraine, Nancy/Metz, France}

    \cortext[cor1]{Corresponding author}

    \begin{abstract}
        The growth of dendritic grains during solidification is often modelled using the Grain Envelope Model (GEM), in which the envelope of the dendrite is an interface tracked by the Phase Field Interface Capturing (PFIC) method. In the PFIC method, an phase-field equation is solved on a fixed mesh to track the position of the envelope. While being versatile and robust, PFIC introduces certain numerical artefacts. In this work, we present an alternative approach for the solution of the GEM that employs a Meshless (sharp) Interface Tracking (MIT) formulation, which uses direct, artefact-free interface tracking. In the MIT, the envelope (interface) is defined as a moving domain boundary and the interface-tracking nodes are boundary nodes for the diffusion problem solved in the domain. To increase the accuracy of the method for the diffusion-controlled moving-boundary problem, an \h-adaptive spatial discretization is used, thus, the node spacing is refined in the vicinity of the envelope. MIT combines a parametric surface reconstruction, a mesh-free discretization of the parametric surfaces and the space enclosed by them, and a high-order approximation of the partial differential operators and of the solute concentration field using radial basis functions augmented with monomials. The proposed method is demonstrated on a two-dimensional \h-adaptive solution of the diffusive growth of dendrite and evaluated by comparing the results to the PFIC approach. It is shown that MIT can reproduce the results calculated with PFIC, that it is convergent and that it can capture more details in the envelope shape than PFIC with a similar spatial discretization.
    \end{abstract}

    %Research highlights
    \begin{highlights}
        \item Novel method for the solution of the Grain Envelope Model (GEM) for dendritic solidification employing sharp interface tracking built on meshless principles.

        \item The new method uses \h-adaptive spatial discretization for differential operator approximation and solute concentration field approximation on a moving boundary problem.

        \item Meshless interface tracking (MIT) 
        %is numerically more accurate than the interface capturing used in previous work because it 
        \hl{uses a conceptually different approach to determine the position of the moving interface than the interface capturing used in previous work. It thus avoids numerical artefacts linked to the use of a phase field for interface capturing.}

        \item Compared to the interface capturing method, MIT requires less computational nodes to capture the details in the dendritic envelope shape.

    \end{highlights}


\end{frontmatter}

% \linenumbers
%
% Introduction
\section{Introduction}
\label{sec:introduction}

Dendritic grains are a prevalent type of crystal growth morphology observed in solidification of metallic alloys. They appear widely in solidification processing: additive manufacturing (3D printing), casting, welding, brazing, etc.\ The size and shape of the grains are a key factor for the properties of the resulting material. Dendritic solidification is controlled by the extraction of latent heat, by diffusion of solute into the liquid surrounding the grains, and by properties of the solid-liquid interface. Therefore it is a multiscale phenomenon and requires the use of models at different scales, ranging from the scale of the dendrite branch (microscopic, typically $10^{-6}\,$m), over the scale of an ensemble of grains (mesoscopic, typically $10^{-3}\,$m) to the process (macroscopic, typically $1\,$m). 

At the scale of an ensemble of dendritic grains, the mesoscopic Grain Envelope Model (GEM)~\cite{steinbach1999three,delaleau2010mesoscopic,souhar2016three} is established as a powerful simulation tool. The GEM describes a dendritic grain by its envelope, which is a smooth surface that circumscribes the tips of the growing dendrite branches within a relatively simple shape. The GEM has been used to predict and characterize shapes of equiaxed dendrites~\cite{steinbach1999three,souhar2016three}, solutal interactions between dendrites~\cite{steinbach2005transient,delaleau2010mesoscopic,olmedilla2019quantitative,chirouf2023investigation}, formation of columnar dendritic patterns~\cite{viardin2017mesoscopic}, and interactions with convection of liquid~\cite{viardin2020mesoscopic}. Furthermore it has been used for upscaling to macroscopic models: macroscopic constitutive laws of grain growth kinetics have been formulated from GEM simulations of the growth of an ensemble of grains~\cite{torabirad2019upscaling}. The capabilities of the GEM to perform quantitative predictions have been assessed extensively by comparisons to experiments and microscopic models~\cite{steinbach1999three,steinbach2005transient,delaleau2010mesoscopic,souhar2016three,olmedilla2019quantitative,viardin2017mesoscopic,viardin2020mesoscopic,tourret2020comparing,boukellal2023multiscale}.

%Dendritic grains are a prevalent type of growth morphology observed in the solidification of metallic alloys. Their complex tree-like appearance results from the growth of individual dendrites, which is primarily governed by the interactions between adjacent grains. The growth of dendritic grains is often modelled using the so-called Grain Envelope Model (GEM)~\cite{souhar2016three,steinbach2005transient} -- a model capable of simulating a large number of grains as a whole, including their collective interactions. The governing concept underlying the GEM is that it is not necessary to describe each branch of the dendrite individually. Rather, it is sufficient to represent the dendrite by its envelope, which is a smooth surface enclosing an entire group of growing dendrite tips within a relatively simple shape. The growth rate of the envelope is determined using a closed-form stagnant-film model that establishes a correlation between the growth rate and the solute concentration near the envelope. Such formulation allows us to use a much coarser numerical discretization (by an order of magnitude coarser) compared to methods describing the solid-liquid interface (e.g., the phase field).

The accuracy of numerical solutions of the GEM~\cite{souhar2016three,tourret2020comparing} is determined primarily by the resolution of the solute concentration in the vicinity of the evolving grain envelope, as well as by the numerical method used to track the envelope. The standard formulations of the GEM ~\cite{steinbach1999three,steinbach2005transient,souhar2016three} use a phase-field-like interface capturing (PFIC) method~\cite{sun2007sharp} to track the grain envelope on a fixed mesh, where the tracked envelope is given by the level set of a continuous phase indicator field. The main advantage of such an approach is that it avoids explicit tracking of the envelope and instead propagates the phase field by solving a transport equation on a \emph{fixed domain}. This results in a relatively straightforward and computationally efficient method that has also been shown to be versatile and robust~\cite{sun2007sharp}.

There are, however, some drawbacks to using PFIC. In the PFIC approach, the envelope is determined from the continuous field that is propagated by solving the transport equation. Therefore, the propagation speed is not calculated only at the front, but needs to be calculated as a scalar field over the entire domain. 
This is difficult to do with high accuracy. To achieve sufficient accuracy of the propagation speed calculation for the PFIC, Souhar et al.~\cite{souhar2016three} introduced a front reconstruction technique, which adds significant computational cost to the method. A further aspect requiring particular attention is the tendency of the PFIC to smoothen the envelope~\cite{sun2007sharp}. This could in some cases hinder the development of physically significant envelope protrusions, for example such that lead to the formation of new branches of the grain.
%The front speed in GEM depends on the positions of the front and the stagnant-film, therefore, accurate determination of the distances from the envelope to the stagnant-film from any mesh point is crucial to obtain an accurate velocity field. Furthermore, the stabilization term in purely advective phase field transport smoothens the interface~\cite{sun2007sharp} and potentially hinders the development of envelope protrusions that could be physically significant. Finally, the standard fixed mesh approach also limits the resolution of small curvature radii, i.e.\ the accuracy of the shape representation is directly related to the fixed mesh spacing. As a result, small-scale features of the phase front, such as protrusions, cannot be accurately resolved on fixed meshes.

If not handled carefully, these aspects can cause undesired numerical artefacts that are potentially detrimental when describing highly non-linear phenomena such as the formation of new branches of the dendrite envelope or the evolution of branches in the initial stages of dendrite growth. It has been shown that the representation of such branching by the GEM remains faithful to the physics~\cite{viardin2017mesoscopic,viardin2020mesoscopic} but the role of numerical artefacts associated with the envelope tracking method has not yet been verified.
\hl{
A way of improving the accuracy of the front propagation without prohibitively increasing the computational cost is to refine the discretization in the vicinity of the front using adaptive spatial discretization}~\cite{Zhang2019d}. \hl{In modeling of dendritic growth such adaptive methods have been employed for solving microscopic thermodynamics-based phase-field models that describe the evolution of the solid-liquid interface of a dendrite during solidification. They have been used in conjunction with finite-element}~\cite{Provatas1999, Sarkis2016}, \hl{finite-difference}~\cite{Greenwood2018, Sakane2022, Guo2015a, Ham2024}, \hl{and meshless}~\cite{Dobravec2022, ghoneim2020smoothed, bahramifar2022local, ghoneim2016new} \hl{methods. Artefacts linked to the PFIC can also be omitted by using a different front propagation method. For dendritic growth several approaches based on level-set methods have been proposed}~\cite{Gibou2003, Tan2006, Ghoneim2018, Ramanuj2019, Limare2023}, \hl{however they also present certain limitations in the representation of curved surfaces}~\cite{DuChene2008}. \hl{Furthermore, methods based on moving tracker particles over a fixed mesh have been developed}~\cite{Juric1996, Reuther2014a}. \hl{Their advantage is a straightforward and accurate representation of the front, but their complexity lies in the coupling between the front tracking points and the underlying mesh. Note that all these approaches have been developed for microscopic dendrite growth models that describe the solid-liquid interface. Such approaches have not been used for mesoscopic dendrite envelope models, such as GEM.}

In this paper, we explore the possibility to leverage meshless methods~\cite{nguyen_meshless_2008, bayona2019insight, tolstykh2003using} to develop an alternative numerical treatment of the GEM with the aim of reducing the aforementioned numerical artefacts. One of the key properties of meshless methods is that they operate on unstructured scattered nodes~\cite{slak2019generation}, which makes the spatial discretization more flexible and therefore particularly suitable for problems with moving boundaries and problems with high spatial variation in the desired accuracy of the solution~\cite{ortega_meshless_2013}, both of which apply to GEM. 

Instead of using a fixed mesh throughout the simulation and using a phase field to capture the interface, we propose to consider the GEM as a moving boundary problem. In this approach the computational domain is the liquid surrounding the grain and the envelope is a moving boundary, where an appropriate boundary condition for the solution of solute diffusion in the liquid is applied. The envelope is described by boundary nodes that are used both as discretization nodes and as trackers of the envelope motion. In addition, we propose to use an adaptive node refinement~\cite{slak2019adaptive} that discretizes the solute field on and near the envelope with a much finer node distribution than in the regions far from the envelope. In summary, we propose an alternative numerical solution method for the GEM using an \h-adaptive meshless interface tracking (MIT) solution procedure.

Since MIT is an interface tracking method, it avoids the complexity of having to provide a velocity field over the entire domain, as well as computing an additional phase-field equation. Instead, the envelope speed is calculated only on the front, more precisely in the nodes that discretize the envelope. Additionally, the meshless methods --- by their design --- also provide all the tools necessary to accurately approximate the solute concentration required to compute the envelope velocity. MIT also avoids the smoothing of protrusions, since the envelope is described directly by trackers, thus omitting artefacts induced by a phase field. Finally, the adaptive spatial discretization offers the possibility to describe small radii of curvature on the envelope and also improves the accuracy of the concentration field in the vicinity of the envelope.

We first introduce the GEM (Section~\ref{sec:model}) and we then discuss its standard PFIC solution along with its limitations (Section~\ref{sec:diffuse}). In the next step (Section~\ref{sec:sharp}), we describe the meshless concept and the proposed MIT solution procedure, where we also discuss the expected numerical error for each principal component of the solution procedure separately. In Section~\ref{sec:istropy}, the entire solution procedure is verified on a synthetic case of isotropic envelope growth. Finally, in Section~\ref{sec:results}, we present a GEM simulation of dendritic growth using MIT and we compare the results with PFIC solutions.

%The use of sharp-interface tracking combined with an adaptive spatial discretization would be advantageous for two reasons: (1) a sharp-interface description would enable us to eliminate the numerical artifacts related to a diffuse-interface; (2) adaptive methods can automatically refine the spatial discretization where required by the form of the solution, typically in the regions of high gradients surrounding the grains, thus increasing the accuracy of the numerical solution.

%
% GEM
\section{Mesoscopic Grain Envelope Model (GEM) for equiaxed isothermal solidification}
\label{sec:model}

The GEM represents a dendritic grain by its envelope. The dendrite envelope is an artificial smooth surface that connects the tips of the actively growing dendrite branches. The idea behind this simplified representation of a dendrite is that the solute interactions between the grains can be accurately simulated without the need for a detailed representation of the branched structure of the solid-liquid interface. Because of the simplified shape, the computational cost of the GEM is several orders of magnitude smaller than that of models that represent the branched structure of the dendrite in detail~\cite{tourret2020comparing,viardin2020mesoscopic}.

\begin{figure}[h]
    \centering
    \includegraphics[width=0.6\textwidth]{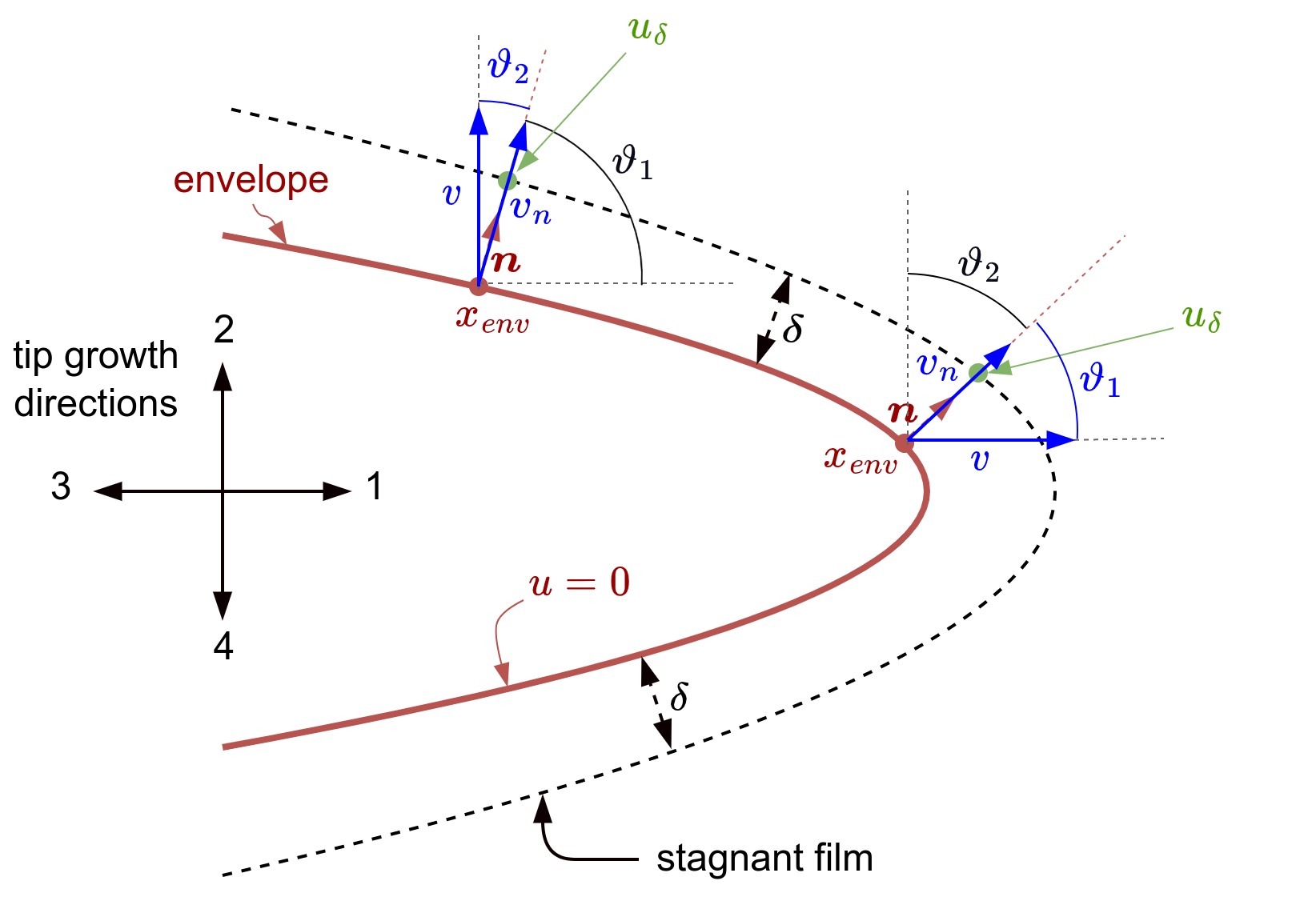}
    \caption{Schematic illustration of the main concepts of the GEM.}
    \label{fig:gem}
\end{figure}

Figure~\ref{fig:gem} illustrates the key concepts of the GEM. The grain is delimited by the dendrite envelope and its growth is controlled by solute diffusion from the envelope into the adjacent, fully liquid domain. The grain contains both liquid and solid phases. The liquid within the grain and on the envelope is assumed to be in a state of thermodynamic equilibrium. In the case of a binary alloy, this means that the solute concentration is determined by the temperature. In the isothermal system considered in this paper, the temperature is uniform and constant throughout the entire domain (i.e., over several dendrites), resulting in a uniform and constant solute concentration on the envelope. It should be noted that the dimensionless model formulation presented below is only applicable to isothermal solidification, as it assumes a uniform and constant solute concentration on the envelope. More general model descriptions can be found in the literature~\cite{delaleau2010mesoscopic,viardin2020mesoscopic}.

The GEM assumes that dendrite branches grow in predefined growth directions\footnote{For example, in $\langle 1 0 0 \rangle$ crystal directions for a typical cubic crystal.}. The growth speed $v_n$ of the dendrite envelope is then given by the speed of the tips $v$ of the dendrite branches by the relation
\begin{equation}
    \label{eq:envelopespeed}
    v_n = v \cos(\vartheta),
\end{equation}
where $\vartheta$ is the smallest of the angles between the outward normal $\b n$ to the envelope and the individual growth directions of the tips. In the present study, the branches are given $\left\langle 10 \right\rangle$ directions, i.e., four possible directions within a two-dimensional domain. The tip speed is calculated from a stagnant-film formulation of the 2D Ivantsov solution~\cite{cantor1977dendritic} which relates the Péclet number $\text{Pe}$ of the tip to the normalized dimensionless concentration $u_\delta$ of the liquid at a distance $\delta$ from the tip, as
\begin{equation}
    \label{eq:cantorvogel}
    u_\delta = \sqrt{\pi \text{Pe}}  \exp(Pe)  \left[ \operatorname{erfc} \left( \sqrt{\text{Pe}} \right) - \operatorname{erfc} \left( \sqrt{\text{Pe}+\delta} \right) \right]
\end{equation}
and from the so-called tip selection criterion that in dimensionless form reads
\begin{equation}
    \label{eq:tipspeed}
    v = (\text{Pe}/\text{Pe}_\text{Iv})^2,
\end{equation}
where $\text{Pe}_\text{Iv}$ is a constant\footnote{$\text{Pe}_\text{Iv}$ is the Péclet number of a free dendrite tip, i.e., a tip growing into an infinite liquid with far-field concentration $u_0$, given by the solution of $u_0 = \sqrt{\pi \text{Pe}_\text{Iv}} \exp(\text{Pe}_\text{Iv}) \operatorname{erfc} \left( \sqrt{\text{Pe}_\text{Iv}} \right)$. Here, solute concentration $u_0$ is a physical parameter of the model.}. Equations~\eqref{eq:cantorvogel} and~\eqref{eq:tipspeed} are used to solve for $v$ at any given point on the envelope. The concentration $u_\delta$ is given by the concentration field in the liquid around the envelope, which is resolved numerically from the diffusion equation
\begin{equation}
    \label{eq:diffusion}
    \frac{ \partial u}{\partial t} = \lap u .
\end{equation}
Equations~(\ref{eq:envelopespeed}--\ref{eq:diffusion}) define the GEM, where the stagnant-film thickness, $\delta$, is a model parameter and $\text{Pe}_\text{Iv}$ is a physical parameter, function of the far-field concentration, $u_0$. The boundary condition for Equation~\eqref{eq:diffusion} on the envelope is $u = 0$. The model also considers that the liquid concentration inside the envelope (grain) is $u=0$.
\hl{
    This means that where the stagnant film overlaps with a neighbouring grain, the stagnant-film concentration, $u_\delta$, becomes $0$ and then the velocity given by Equation
}
~\eqref{eq:tipspeed}
\hl{
    is also zero. The envelope growth therefore stops when two dendrite envelopes grow near to a distance equal to the stagnant-film thickness, $\delta$. In order for the envelope velocity model to provide a physically correct growth velocity, $\delta$ needs to be of a certain thickness and should not be significantly smaller or larger (see refs.
}
~\cite{souhar2016three}, \cite{tourret2020comparing}).
\hl{
    In the present work the dimensionless thickness of the stagnant film is $\delta=1$.
}

In the simulations presented in this paper the solute diffusion flux across the outer domain boundary is zero, i.e., the boundary condition for Equation~\eqref{eq:diffusion} on the outer liquid boundary is $\b n \cdot \nabla u = 0$, where $\b n$ is the normal vector to the outer boundary $\Gamma_\ell$. The initial concentration of the liquid is $u = u_0$, i.e., the liquid is homogeneous and has a dimensionless supersaturation of $u_0$. In this work we use a stagnant-film thickness of $\delta = 1$.

Note that all equations are written in terms of dimensionless quantities, where the dimensionless concentration is defined by $$u=\frac{1-\frac{C}{C^\star}}{1-k},$$ with the dimensional concentration $C$ , the equilibrium concentration $C^*$ of the liquid at the given solidification temperature and the equilibrium solid-liquid solute partition coefficient $k$. The other characteristic scales for normalization are $v_\text{Iv}$ for velocity, $D/v_\text{Iv}$ for length, $D / v_\text{Iv}^2$ for time, where $v_\text{Iv}=4\sigma^* D \text{Pe}_\text{Iv}^2/d_0$ is the steady-state velocity of the free tip, $\sigma^*$ is the tip selection parameter, $D$ is the diffusion coefficient, and $d_0$ is the capillary length.

%
% Diffuse-interface
\section{GEM with phase-field interface capturing on a fixed mesh (PFIC)}
\label{sec:diffuse}

The standard formulations of GEM~\cite{souhar2016three,steinbach2005transient,steinbach1999three} use a PFIC method~\cite{sun2007sharp} to track the grain envelope on a fixed mesh. In PFIC, the tracked front is given by the level set $\phi=\phi\env$ of a continuous indicator field $\phi$. The transition of $\phi$ between 1 (grain) and 0 (liquid) is smooth but compact; it follows a hyperbolic tangent profile with characteristic width $W_\phi$~\cite{sun2008atwo,souhar2016three}
\begin{equation}
    \label{eq:diffuse-kernel}
    \phi(n) = \frac{1}{2} \left[ 1 - \tanh \left( \frac{n}{2 W_\phi} \right) \right] \,,
\end{equation}
where $n$ is the distance from the center of the hyperbolic tangent profile. The evolution of $\phi$ is given by a phase-field equation that ensures the transition is self-preserving and retains its shape and width
\begin{equation}
    \label{eq:diffuse}
    \frac{\partial \phi}{\partial t}
    + v_n \b n \cdot \nabla \phi =
    -b
    \underbrace{
        \left[
            \lap \phi
            - \frac{\phi(1 - \phi)(1-2\phi)}{W_\phi^2}
            - \left| \nabla \phi \right| \nabla \cdot \left( \frac{\nabla \phi}{\left| \nabla \phi \right|} \right)
            \right]
    }_{\mathrm{stabilization \ term}},
\end{equation}
where $v_n$ is the envelope growth speed defined as a scalar field. The term on the right hand side of the equation is a stabilization term that ensures the phase-field retains the hyperbolic tangent transition. The coefficient $b$ is a numerical parameter that controls the relaxation of the phase-field profile. Generic criteria on how to choose $b$ and $W_\phi$ in order to minimize the error and to ensure the stability of the PFIC method are discussed by Sun and Beckermann~\cite{sun2007sharp}. Specific guidelines for the use of the PFIC for envelope tracking in the GEM are provided by Souhar et al.~\cite{souhar2016three}. In this work the parameters $\phi\env$, $W_\phi$ and $b$, as well as the grid size and time step are selected following the guidelines of Souhar et al.~\cite{souhar2016three}. Note that the choice of $\phi\env=0.95$ is to favour the accuracy of the tip of the envelope in a convex shape.

%FROM sun2008atwo:
%The term in the square brackets on the right-hand side of Equation (24) is added to the continuity equation, Equation (15), to sustain the hyperbolic tangent φ profile across the diffuse interface during motion. The last term in the square brackets is the so-called counter-term [4] that cancels out curvature-driven interface motion at leading order for a finite interface width [see also Equation (6)]. For a flat interface where κ = 0, the last term in the square brackets is equal to zero and the first two terms (i.e. ∇2
%φ−φ (1 − φ) (1 − 2φ) /δ2) yield the hyperbolic
%tangent φ profile as the stationary solution of Equation (24). During interface motion, the term in the square brackets sustains the hyperbolic tangent φ profile, while the ∇2
%φ term smoothes out
%interface singularities.

The main advantage of this computational method is that it avoids explicit tracking of the envelope and instead solves the PDE from Equation~\eqref{eq:diffuse} on a fixed mesh. Additionally, Sun and Beckermann~\cite{sun2007sharp} demonstrated its versatility and robustness as a general front capturing method. When coupled with the GEM, however, a few aspects of the method require further attention.
\begin{description}
    \item[Accurate calculation of the front speed.]
        To propagate the phase field with Equation~\eqref{eq:diffuse}, the front speed, $v_n$, needs to be provided everywhere where $\nabla \phi$ is non-negligible, i.e., in the vicinity of the front.
        %Because the front speed depends on the positions of the front and of the stagnant-film, accurate determination of the distances to the envelope and to the stagnant-film from any mesh point $\x \in \Omega$ is crucial to obtain an accurate speed field.
        In the PFIC-based GEM the speed $v_n$ at point $\b x\in \Omega$ is equal to the speed of the closest point lying on the sharp front, i.e., $v_n (\b x) = v_n (\b x\env)$.
        %Approximations:
        %(1) CLOSEST MARKER POINT = CLOSEST POINT
        %(2) NORMALIZED GRADIENT IN x = NORMAL AT THE FRONT
        %(3) INTERPOLATION of u in x\delta by a second-order method
        Note that by definition the vector between a point and the closest point on the front, $\b x - \b x\env$, is normal to the front in $\b x\env$.
        The speed of point $\b x\env$ depends on the concentration at a given distance $\delta$ from the front, $u_\delta =u(\x_\delta)=u(\b x\env + \delta \, \b n)$, where $\b n$ is the normal vector to the front at point $\b x\env$ and $\x_\delta$ is the stagnant-film position. Thus, the speed of point $\b x\env$ is given by Equations~\eqref{eq:envelopespeed}, \eqref{eq:cantorvogel}, and~\eqref{eq:tipspeed} as a function of $u_\delta$, $\delta$, and $\b n$. Because the speed in $\b x$ is equal to the speed in $\b x\env$, we can write a functional expression $v_n(\b x) = f(\b x - \b x\env, u(\b x\env + \delta \, \b n), \delta)$, which shows the three quantities that need to be approximated in order to compute the frontal speed:
        \begin{enumerate}[(i)]
            \item the distance $\b x - \b x\env$ from the front,
            \item the vector $\b n$ normal to the front in $\b x\env$ and
            \item the interpolation of the concentration $u_\delta$ to the corresponding point on the stagnant-film $\b x_\delta =\b x\env + \delta \, \b n$.
        \end{enumerate}

        The accuracy of these three approximations determines the accuracy of the calculated speed, $v_n$.
        % It is important to note that $f$ is a strictly increasing function of $u_\delta$ for any given $\delta$. Also note that in the vicinity of the envelope the dimensionless concentration, $u$, strictly increases with distance from the envelope. This means that an error that overestimates the distance of the stagnant-film from the envelope, results in an overestimated speed of the envelope, and vice versa.
        In this work the PFIC-based GEM uses the front reconstruction method proposed by Souhar et al.~\cite{souhar2016three} to determine the distance $\b x - \b x\env$ from the front and the normal vector, $\b n$. In this method, the envelope is reconstructed by marker points densely distributed over the surface defined by the level set $\phi=\phi\env$. The distance from a point $\b x$ to the envelope is then approximated by the projection of the distance to the closest marker, $\b x_\text{m}$, on an approximated normal vector $$|\b x - \b x\env| \approx - (\b x - \b x_\mathrm{m}) \cdot \b n,$$ where the normal is approximated by a normalized gradient of the phase field in $\b x$ with the following expression
        \begin{equation}
            \label{eq:grad_normal}
            \b n \approx \frac{\nabla \phi(\b x)}{ \lvert \nabla \phi(\b x) \rvert}.
        \end{equation}
        The same normal vector approximation is used to calculate the stagnant-film point, $\b x_\delta$. The interpolation of $u_\delta$ in $\b x_\delta$ is further done by a second-order interpolation scheme.

        Note that, in principle, the distances can be determined via the phase field, $\phi$. However, as pointed out by Souhar et al.~\cite{souhar2016three}, such approach results in a self-reinforcing error and is not sufficiently accurate for practical use.
        %The front reconstruction method provides much better accuracy. Its disadvantage is that it increases the computation time considerably because the front reconstruction with markers needs to be done at each timestep.
        %
    \item[Damping of protrusions.]
        The stabilization term of Equation~\eqref{eq:diffuse} smoothens interface singularities due to the dissipative nature of the $\nabla^2 \phi$ term~\cite{sun2007sharp}. The damping of ripples on the envelope is controlled by the relaxation coefficient $b$. Unfortunately, it is not known to what extent such damping affects the development of envelope protrusions that have a physical meaning.
    \item[Resolution of small curvature radii.]
        The phase-field method can accurately describe convex radii of front curvature larger than $W_\phi \left[ 6 - \ln \left( \frac{\phi\env}{1-\phi\env} \right) \right]$ and concave radii larger than $W_\phi \left[ 6 + \ln \left( \frac{\phi\env}{1-\phi\env} \right) \right]$~\cite{souhar2016three}. With standard mesh spacing of $h\leq W_\phi / \sqrt{2}$~\cite{souhar2016three,sun2007sharp}, this means that the resolution of shape representation is directly proportional to the fixed mesh spacing, implying that small-scale front features like protrusions may not be adequately resolved on fixed meshes.
        %Also, the radius of the initial envelopes must be sufficiently large if one wants to exclude any errors linked to the front cacpturing method. 
        %
\end{description}

\section{GEM with meshless interface tracking (MIT) on \h-refined scattered nodes}
\label{sec:sharp}
Using the meshless principle for the solution of the GEM opens a variety of possibilities for spatial discretization.
%MZ One is to maintain the grain envelope populated with nodes ($\Gamma_\mathrm{e}$) throughout the simulation, which effectively means that instead of reconstructing the envelope from evolving continuous indicator field, the envelope is described by a set of nodes, in our case boundary nodes, as the envelope defines the internal part of the domain (Figure~\ref{fig:domain_discretization}). We call this a sharp-interface (SI) approach.   
We propose a method that considers the grain envelope as a boundary ($\Gamma_\mathrm{e}$) of the liquid domain ($\Omega$), as schematically shown in Figure~\ref{fig:domain_discretization}. The method adapts the node distribution to ensure that the envelope boundary is densely populated with nodes throughout the simulation. These nodes are then also used for the tracking of the envelope. This effectively means that the grain envelope is described as a sharp interface by a set of discretization nodes -- in this case, the boundary nodes from $\Gamma_\mathrm{e} \subset \partial \Omega$. The internal part of the domain is then the liquid surrounding the grain, i.e., the area between the envelope and the outer liquid boundary, $\Gamma_\ell \subset \partial \Omega$ (see Figure~\ref{fig:domain_discretization} for clarity). Equation~\eqref{eq:diffusion} is solved in the domain.
%We refer to this approach as the sharp-interface (SI) approach.
%level set indicator field, which we can avoid using meshless method.  i.e. the envelope is not reconstructed from the evolving continuous indicator, but it is the boundary of the domain.    
%In the context of sharp-interface tracking methods, the dendrite envelope is defined at any time during the simulation by a set of discretization nodes $\Gamma_\mathrm{e}$ positioned at the solidification front as shown in FIgure~\ref{fig:domain_discretization}. 
%It is our assumption that the numerical artefacts appearing in the diffuse-interface capturing methods, are completely avoided when the envelope shape is sought. 
We use scattered nodes to discretize the domain $\Omega$, which allows us to refine the local field description in the vicinity of the grain envelope and consequently improve its discretization quality. Moreover, since the interpolation of scattered data is the backbone of meshless methods, the sampling of the concentration $u_\delta$ (needed to determine the envelope velocity) is inherently present in the solution procedure. This ensures that the error of the sampling is in the worst case of the same order of magnitude as the error of the discretization of the partial differential operators used in solving Equation~(\ref{eq:diffusion})~\cite{le2023guidelines,tominec2021least}. On a theoretical level, the MIT approach should allow us to mitigate some of the problems that are inherently present in PFIC approach (see Section~\ref{sec:diffuse}), but at the cost of a higher computational complexity per node, associated with the use of meshless methods~\cite{jancic_monomial_2021}.

The proposed MIT solution procedure for the GEM starts by building the initial domain -- the initial grain envelope, $\Gamma_\mathrm{e}$, as a circle with radius $r_0$ and the outer liquid boundary, $\Gamma_\ell$ as a square box with side length $L$. The domain and the boundaries are populated with nodes and the initial solute concentration is set to $u_0$ in $\Omega$ and $\Gamma_\ell$ (liquid) and to $0$ in $\Gamma_\mathrm{e}$ (envelope). Afterwards, the explicit forward time marching starts to simulate the solute diffusion and the spatial expansion of the envelope. The entire MIT solution procedure is given in Algorithm~\ref{alg:simulate} and its elements are explained in more detail in the following subsections.

\begin{algorithm}[h]
    \caption{MIT solution procedure for simulation of dendritic grain growth.}
    \label{alg:simulate}
    \vspace{1mm}
    \textbf{Input:} GEM, nodal density function $h:~\Omega\to \R$, RBF-FD approximation basis $\xi$, boundary concentration $u_0$, number of time steps $I_{\text{max}}$, list of cumulative boundary node displacements $\varphi$, time step $\mathrm{d}t$, parametric spline degree $k$.\\
    \textbf{Output:} Concentration field $u$ and domain discretization $\Omega$. \\

    \begin{algorithmic}[1]
        \State $\Omega, u \gets \Call{initialize}{h, u_0}$
        \Comment{Obtain initial domain $\Omega$ and initial concentration field $u$.}
        \State $\Gamma_{e} \gets \Call{get\_envelope}{\Omega}$
        \Comment{Obtain envelope boundary nodes.}
        \State $\varphi \gets \b 0$
        \Comment{Initialize cumulative displacements of grain envelope nodes.}\\

        \Function{simulate\_growth}{GEM, $\Omega, h, \xi, \varphi, u, \mathrm{d}t$}
        \For{$i$}{0}{$I_{\text{max}}$} \label{alg:while}
        \State $\L \gets \Call{operator\_approximation}{\text{GEM}, \Omega, \xi}$
        \Comment{Approximation of differential operators.}
        \State $u \gets \Call{solve\_diffusion}{\Omega, \L, \mathrm{d}t}$\label{alg:line_diffusion}
        \Comment{Solute diffusion problem.}
        \State $\Gamma_\mathrm{e}^\star \gets \Call{new\_envelope}{\text{GEM}, \Omega, \Gamma_\mathrm{e}, u,\xi, \mathrm{d}t}$ \label{alg:new_envelope}
        \Comment{Obtain new envelope nodes positions as defined by GEM.}
        \\
        \If{$\Call{max\_tip}{\Gamma_\mathrm{e}^\star} >  \Call{max\_tip\_allowed}{}$}
        \State \Return $u, \Omega$
        \Comment{End simulation.}
        \EndIf
        \\
        \State $\varphi \gets \Call{add\_displacements}{\Gamma_\mathrm{e}, \Gamma_\mathrm{e}^\star, \varphi}$
        \Comment{Obtain list of updated total displacement magnitudes.}

        \If{$\Call{max}{\varphi} > \he / 2$}\label{alg:disp}
        \State $\Omega^\star \gets \Call{re-discretize\_domain}{\Gamma_\mathrm{e}^\star, h, k}$\label{alg:re-discretize}
        \Comment{Construct new domain $\Omega$ using envelope nodes $\Gamma_\mathrm{e}^\star$.}
        \State $u \gets \Call{map\_concentration}{\Gamma, \Omega^\star, u, \xi}$ \label{alg:map}
        \Comment{Map concentration field to new discretization $\Omega^\star$.}
        \State $\Gamma_{e} \gets \Call{envelope\_boundary}{\Omega}$
        \Comment{Obtain envelope boundary nodes.}
        \State $\Omega \gets \Omega^\star$
        \Comment{Re-assign domain.}
        \State $\varphi \gets \b 0$
        \Comment{Reset cumulative displacements of grain envelope nodes.}
        \Else
        \State $\Omega, u \gets \Call{remove\_nodes}{h, \Omega}$
        \Comment{Remove nodes from $\Omega$ and $u$ closer to envelope than $(0.8 \he)^2$.}
        \EndIf

        \EndFor \\
        \State \Return $u, \Omega$
        \EndFunction
    \end{algorithmic}
\end{algorithm}
\subsection{Domain discretization}
\label{sec:discretization}
% \mz{Tale odstavek spodaj sem malo poenostavil. Izvirno verzijo sem bolj težko razumel, zato sem preformuliral in tudi izpustil nekaj podrobnosti, ki se mi ne zdijo bistvene. Detajli o algoritmu so itak v paperju Slak\&Kosec, 2019.}
To discretize the domain $\Omega$ with scattered nodes, we employ a dedicated dimension-independent variable density meshless node positioning (DIVG) algorithm~\cite{slak2019generation}. DIVG is an iterative algorithm that populates the entire domain in an advancing front manner, starting from one or more initial seed nodes stored in the \emph{expansion queue}. By default, the algorithm uses all boundary nodes, i.e.\ nodes from $\Gamma_\mathrm{e}$ and $\Gamma_\ell$, as seed nodes. In each iteration, a single node $\b{x}_i$ is removed from the queue and is \emph{expanded}. Expansion means that candidates for new nodes are uniformly positioned on
%the annulus with inner radius $h_i(\b{x}_i)$ and outer radius $h_o(\b{x}_i)$ centred at $\b{x}_i$. Note that the difference between inner and outer radius is relatively small allowing to safely assume $h_i(\b{x}_i)=h_o(\b{x}_i)=h(\b{x}_i)$ for further discussion. 
a circle with radius $h(\b{x}_i)$, then those that fall outside the domain (we discuss the \texttt{inside check} in Section~\ref{sec:reconstruction}) or lie too close to the already accepted nodes, are discarded. The remaining candidates are accepted, i.e., they are added to the domain and to the expansion queue. The iterative procedure continues until no new nodes can be added to the domain and the expansion queue is empty.

\begin{figure}
    \centering
    \includegraphics[width=0.5\textwidth]{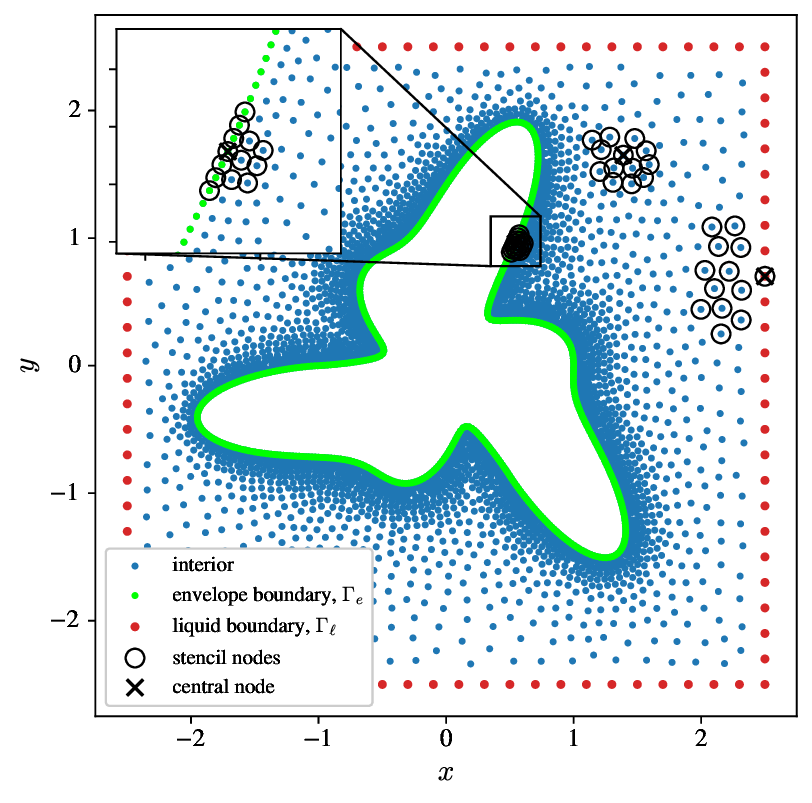}
    \caption{Example \h-refined domain discretization. Additionally, example stencils on Dirichlet (green), Neumann (red) boundary and in the interior (blue) are shown.}
    \label{fig:domain_discretization}
\end{figure}

A detailed description of the DIVG algorithm can be found in~\cite{slak2019generation}. The key feature of DIVG is its ability to populate the domain with a spatially variable node spacing, $h(\b{x})$, and that it can populate complex geometries regardless of the dimensionality of the domain~\cite{depolli_parallel_2022, de2019fast}. An example of nodes generated with DIVG within $\Omega$ and for its Dirichlet (envelope $\Gamma_\mathrm{e}$) and Neumann (liquid $\Gamma_\ell$) boundaries is shown in Figure~\ref{fig:domain_discretization}.

To ensure an accurate description of the grain envelope and of the steep diffusion field in the vicinity of the envelope we employ \h-refinement through the following internodal distance function
\begin{equation}
    \label{eq:slakNodes}
    h(\x) = \begin{cases}
        \he                                                      & , \text{ if } \x \in \Gamma_\mathrm{e} \\
        h_\ell                                                   & , \text{ if } \x \in \Gamma_\ell       \\
        \he + (h_\ell - \he)  \frac{ d(\x)}{ d_\mathrm{closest}} & , \text{ otherwise}
    \end{cases},
\end{equation}
where $\he$ and $h_\ell$ are the internodal distances on the grain envelope and the outer liquid boundary, $d(\x)$ is the Euclidean distance between $\x$ and the nearest node on the envelope, and $d_\text{closest}$ is the Euclidean distance between the node on the envelope boundary, $\Gamma_\mathrm{e}$, nearest to $\x$ and the node on the outer liquid boundary, $\Gamma_\ell$, nearest to $\x$. In our case, $\he < h_\ell$. Such definition of $h(\x)$ gives a linearly decreasing internodal distance towards the envelope, as shown in Figure~\ref{fig:domain_discretization} with green coloured nodes on an illustrative dendrite-like shape.

%
%
%\subsubsection*{Criterion for domain re-discretization}
\subsection{Criterion for domain re-discretization}
With the growth of the grain the envelope boundary moves and the shape and size of the domain therefore evolves. The spatial discretization must therefore be adapted with time. Besides the approximation of the differential operators, the domain discretization is among the most computationally demanding parts of Algorithm~\ref{alg:simulate}. Since the displacements of the envelope nodes within a single time step are typically small (much smaller than the local internodal distance $\he$) and since RBF-FD is relatively robust to non-optimal nodal positions~\cite{jancic-stability,le2023guidelines}, it is not necessary to regenerate the discretization of the computational domain at each time step. Instead, we use a criterion based on the cumulative displacement magnitudes for each node $\x\env$ on the envelope.
The cumulative displacement magnitudes are collected in an array-like data type $\varphi$, which essentially holds information about the total displacement magnitude of all envelope nodes since the last domain re-discretization.

Based on the maximum cumulative displacement among all nodes, $\max(\varphi)$, we then distinguish between two possible cases (see line~\ref{alg:disp} in Algorithm~\ref{alg:simulate}):
\begin{description}
    \item[$\max(\varphi) < \he /2:$] The largest cumulative displacement of the envelope node is smaller than half of the local internodal distance on the envelope, $\he$. Note that $h_\mathrm{e}$ is a parameter of the method and is constant. In this case re-discretization of the computational domain is not yet necessary. The envelope nodes are repositioned from $\Gamma_\mathrm{e}$ to $\Gamma_\mathrm{e}^\star$, according to the GEM. To avoid that computational nodes from the interior of the domain are too close to or inside\footnote{A detailed description of the decision-making algorithm weather a given node is inside the envelope or not is provided in Section~\ref{sec:reconstruction}.} the recomputed grain envelope (as this may affect the stability of the RBF-FD approximation), all nodes that fall inside the new envelope $\Gamma_\mathrm{e}^\star$ or are closer than $0.8\he$ to the envelope are removed from the domain $\Omega$.

    \item[$\max(\varphi) \geq \he/2:$] The largest cumulative displacement is larger or equal to half of the local internodal distance $\he$. In this case a complete re-discretization of the computational domain is required. The algorithm for re-discretization is given in Algorithm~\ref{alg:reconstruction} in line~\ref{alg:rediscretize} and gives the new computational nodes. The corresponding concentration field is obtained by mapping (interpolating) the field from the old to the new discretization (Algorithm~\ref{alg:simulate}, line~\ref{alg:map}). The interpolation method is discussed in detail in Section~\ref{sec:pu}.
\end{description}

%Note that whether the re-discretization is required or not is based on the internodal distance $\he$ on the dendritic envelope. While this value could be optimized to increase the computational efficiency even further, it must be noted that the larger the threshold value the larger can be the distortions of the envelope shape. Our observations showed that $\he/2$ is a good value to improve the computational efficiency of the solution procedure with minimal re-construction distortions introduced.
%
%

\subsection{Meshless approximation of differential operators}
\label{sec:approx}
The next step is the discretization of partial differential operators $\L$, in our case $\lap$ and $\nabla$, which appear in Equation~\eqref{eq:diffusion} and on the Neumann boundary, respectively. For each computational node $\x_i$ in the domain, a set of nearby nodes $\mathcal N_i = \left \{ \x_j \right \}_{j=1}^n$, commonly called \emph{stencil} or \emph{support} nodes, is defined. While special stencil node selection algorithms showed promising results~\cite{davydov2022improved,9803651}, we use the simplest approach and choose the closest $n$ nodes according to the Euclidean distance, obtained efficiently by constructing a $k$-d tree. Example stencils in the interior of the domain and on the Dirichlet and Neumann boundaries are shown in Figure~\ref{fig:domain_discretization}.

Once the stencil is defined, a linear differential operator $\L$ in node $\x_c$ is approximated over a set of $n$ stencil nodes $\mathcal N_c$ with
\begin{equation}
    \label{eq:ansatz}
    \L \big (u(\x_c)\big ) \approx \sum_{j=1}^nw_ju(\x_j),
\end{equation}
for stencil nodes $\x_j\in \mathcal{N}_c$, function $u$ and weights $\w$. The weights are obtained by solving a linear system, assembled from the chosen set of basis functions.

In this work, we use polyharmonic splines (PHS)~\cite{bayona2017role} augmented with $N_p=\binom{m+d}{m}$ monomials up to order $m$ in the $d=2$ dimensional domain $\Omega$.

This corresponds to a commonly used meshless method, also referred to as the RBF-FD approximation~\cite{tolstykh2003using,flyer2016role,bayona2017role}.

By enforcing the equality of Equation~\eqref{eq:ansatz} for the approximation basis, we can write a linear system
\begin{equation}
    \label{eq:rbf-system-aug}
    \underbrace{\begin{bmatrix}\b \Phi & \b P \\\b P^\T & \b 0\end{bmatrix}}_{\b A} \begin{bmatrix} \b w \\ \b \lambda \end{bmatrix} = \begin{bmatrix} \b\ell_\phi \\ \b\ell_p \end{bmatrix},
\end{equation}
with matrix $\b P \in \R^{n\times N_p}$ of the evaluated monomials, matrix $\b \Phi \in \R^{n\times n}$ of the evaluated RBFs and $\b \ell_\phi$ and $\b \ell_p$ are vectors of values assembled by applying the considered operator $\L$ to the RBFs and monomials respectively. To obtain the weights $\w$, the system is solved and Lagrangian multipliers $\b \lambda$ are discarded.

\hl{Note that the conditionally positive definite RBFs alone do not ensure the positive definiteness of matrices $\b A$. The interpolation problem becomes positive definite Afterwards monomials are added to the approximation basis}~\cite{bayona2017role}\hl{. Another consequence of adding monomials to the approximation basis is the achieved consistency up to a certain order.} Moreover, we follow the recommendations of Bayona~\cite{bayona2019insight} and define the stencil size $n$ as
\begin{equation}
    n = 2 N_p =2 \binom{m + d}{m}
    \label{eq:stencil_size}
\end{equation}
to ensure a stable approximation.

The approximation~\eqref{eq:ansatz} also holds for $\L \coloneqq 1$, effectively allowing us to use the same formulation for the approximation of field values at locations that do not correspond to the discretization points. This property will be used in Section~\ref{sec:pu}, where we seek an accurate interpolation of the concentration field at the stagnant-film positions $\x^\delta$.

Note that the formulation of RBF-FD allows us to control the order of the interpolation through the order of the augmenting monomials~\cite{jancic_monomial_2021, bayona2017role}. In this study, we will experiment with the second ($m=2$) and fourth ($m=4$) orders.
%
%

%\subsection{The moving boundary problem -- modifying the dendrite envelope positions}
\subsection{Propagation of the envelope -- moving the boundary}
\label{sec:growth}
At each time step, the new solute concentration field, $u$, is obtained by discretizing the diffusion Equation~\eqref{eq:diffusion} with the meshless approximation and solving it explicitly using the forward Euler time-marching scheme (see line~\ref{alg:line_diffusion} of Algorithm~\ref{alg:simulate}). With the known concentration field, the envelope velocities are calculated by the GEM and the grain envelope is then propagated (line~\ref{alg:new_envelope} of Algorithm~\ref{alg:simulate}). Recall that the envelope speed depends on the concentration $u_\delta$ at the stagnant-film and on the direction of the envelope normal (see Equations~(\ref{eq:envelopespeed}--\ref{eq:tipspeed})).
%As a result, the GEM model returns the tip speed, $v$, for each envelope node, $\b x\ienv$, and the considered envelope node is moved to a new position, $\x'\ienv$,
Accordingly, each envelope node, $\b x\ienv$, is propagated to a new position, $\x'\ienv$, by
\begin{equation}
    \label{eq:displacement}
    \x'\ienv = \x\ienv + \cos (\vartheta_i) \, v(u_{\delta,i}) \, \b n_i \, \mathrm{d}t \,,
    \quad \forall \x\ienv \in \Gamma_\mathrm{e},
\end{equation}
where the speed $v$ is given by the GEM model (Equations~(\ref{eq:cantorvogel}--\ref{eq:tipspeed})) as a function of $u_{\delta,i}$, which is the concentration at $\b x_{\delta,i} = \x\ienv + \delta \, \b n_i$ and $\vartheta _i$ is the minimum angle between the envelope normal $\b n _i$ and the 4 unit directional vectors, i.e., $\pm \b {\widehat e}_x$ and $\pm\b{\widehat e} _y$.

It follows that the numerical error in the shape of the grain envelope depends on the accuracy of the evaluation of the stagnant-film concentration, $u_\delta$, and of the normal direction, $\b n_i$. These error sources are studied in the following subsections. $u_\delta$ is calculated by interpolation of the concentration field; the accuracy of different interpolation methods is discussed in Section~\ref{sec:pu}. The normal, $\b n$, is calculated from a parametric reconstruction of the envelope with splines. The same reconstruction is also used for the discretization of the envelope and allows us to ensure quasi-uniform internodal distances on the envelope and consequently to generate a robust spatial discretization of the domain. The envelope reconstruction is discussed in Section~\ref{sec:reconstruction} along with the influence of the spline degree, a key parameter of the reconstruction.

\subsubsection{Interpolation of the solute concentration field}
\label{sec:pu}
%In GEM, the envelope velocity $v_n$ is strongly dependant on the solute concentration $u_\delta$ at the stagnant-film around the envelope. 
In the GEM, the stagnant-film concentration, $u_\delta$, is required at any position $\b x^\delta = \x\env +\delta \b n$ within the domain space $\Omega$. Since the concentration is only given at the discretization points $\b x_i\in \Omega$, a suitable interpolation is required. We consider the following three interpolation methods:
\begin{description}
    \item[Shepard's interpolation,] also known as inverse distance weighing, where the interpolant is a weighted average over the neighborhood,

        i.e., the $n$ points closest to $\b x^\delta$, defined as
        \begin{equation}
            \label{eq:shep}
            u(\b x^\delta ) =
            \begin{cases}
                \dfrac{\sum_{i=1}^n\omega_i(\b x^\delta)u(\x_i)}{\sum_{i=1}^n\omega_i(\b x^\delta)} & ,\text{ if } d(\x, \x_i) \neq 0 \text{ for all } i \\
                u (\x_i)                                                                            & ,\text{ if } d(\x, \x_i) = 0 \text{ for some } i
            \end{cases},
        \end{equation}
        for weights $\omega_i(\b x^\delta) = \left \| \b x^\delta - \x_i \right \| ^{-\alpha}$ with power parameter $\alpha = 2$.
    \item[The RBF-FD interpolation,] which uses the RBF-FD approximation described in Section~\ref{sec:approx} to build a local interpolant $F_i(\x^\delta)$ over stencil nodes $\mathcal N (\x^\delta)$ to obtain a local interpolation of the concentration in $\b x^\delta$ for a chosen approximation basis.
    \item [The partition-of-unity (PU) interpolation,] which uses local RBF-FD interpolants and additionally constructs a global representation of the concentration field by employing the PU interpolation~\cite{pu-interpolation,DEMARCHI2019331}. Essentially, the PU approach builds local interpolants $F_i$ over stencil nodes $\mathcal N_i$ for all $N$ discretization nodes $\x_i \in \Omega$. The global interpolant is then constructed as
          \begin{equation}
              F(\x) = \sum_{i=1}^Nw_i(\x)F_i(\x),
          \end{equation} using compactly supported weights $w_i$ that form a partition-of-unity, $\sum_{i=1}^Nw_i(\x)=1$, for all $\x \in \Omega$. The most efficient way to construct the weights is to use the Shepard's construction (see Equation~\ref{eq:shep}) thus
          \begin{equation}
              w_i(\x) = \frac{\omega(\left \| \x - \x_i \right \| /r_i)}{\sum_{j=1}^N\omega(\left \| \x - \x_j \right \| / r_j)},
          \end{equation}
          where $r$ determines the effective radius of the individual interpolants $F$ -- the larger the radius, the more pronounced the averaging effect. Such a construction effectively leads to a weighted average of the local interpolants $F_i$, where $\left \| \x_i - \x \right \| < r_i$, otherwise the weight is 0.

\end{description}
% \subsubsection*{Performance analysis of the proposed interpolation approaches}
To evaluate the performance of the proposed interpolation methods, all three are further discussed and compared in terms of computational complexity and interpolation accuracy.

\begin{itemize}
    \item \emph{Computational complexity} \newline
          The computational complexity of the three interpolation methods is anything but equivalent. Common to all is the construction of a $k$-d tree that allows us to query the nearest $n$ neighbours. Constructing such a structure for a domain with $N$ nodes requires $\mathcal{O}(N\log N)$ operations, while searching for $n$ nearest neighbours requires $\mathcal{O}(n\log N)$ operations.

          The evaluation of Shepard's interpolant then requires additional $\mathcal{O}(n)$ operations. The construction of the RBF-FD interpolant requires $\mathcal{O}(n^3)$ operations to solve the dense linear interpolation system from Equation~\eqref{eq:rbf-system-aug} and additional $\mathcal{O}(n)$ operations to evaluate. The computational complexity of PU is generally higher than that of RBF-FD. In the worst case, i.e., with uniformly distributed nodes and a single query point, the construction of the PU interpolant requires $\mathcal{O}(n(n^3 + n))$ operations. The first term comes from the construction of the RBF-FD interpolant and the second term comes from the calculation of Shepard's weights. To evaluate the PU interpolant, further $\mathcal{O}(n)$ operations are required.

          Thus, the construction of the PU interpolant is the most computationally expensive, followed by the RBF-FD approximant and finally Shepard's interpolant as the cheapest of the three.
    \item \emph{Interpolation accuracy} \newline
          The accuracy of each interpolation approach is evaluated using a two-dimensional Sibson's function
          \begin{equation}
              f(\b x) = f(x,y) =\cos(4\pi \sqrt{(x-1/4)^2 + (y -1/4)^2})
          \end{equation}
          for the domain space $\b x \in [0,1]^2$. We generate $N_{\mathrm{fit}}$ uniformly scattered fit points on which we construct the three proposed interpolation methods and we evaluate the accuracy of the interpolations on $N_\mathrm{test}\approx10^{5}$ regularly positioned test nodes, ${\b x}' \in [0, 1]^2$. The performance of the three interpolation methods is evaluated in Figure~\ref{fig:interpolation}. On the left, we show the interpolated values using $N_{\mathrm{fit}}=384$ fitting points (marked with black crosses).

          In the center plot, we evaluate the relative error of all three interpolation methods as function of the neighbourhood/stencil size, $n$.
          For clarity, black dashed lines indicate the stencil sizes $n=12$, recommended by RBF-FD for a second-order method ($m=2$), and $n=30$, recommended for a fourth-order method ($m=4$) (see Equation~\eqref{eq:stencil_size}). We find that Shepard's interpolation leads to interpolation errors that are almost an order of magnitude larger than the errors of RBF-FD and PU. Increasing the neighbourhood size, $n$, only worsens the accuracy of Shepard's interpolation. \hl{RBF-FD and PU show comparable accuracy, both at second and fourth order, with the error decreasing with increasing monomial order}.

          Similar observations are made by a convergence analysis with respect to $N_{\mathrm{fit}}$ and with recommended stencil sizes $n$, shown in the right plot of Figure~\ref{fig:interpolation}. The number of closest neighbours for Sheppard's interpolation is set to $n=3$, as it yielded the highest accuracy in the centre plot. The black dashed lines and their slopes, $k$, show the observed order of convergence.

          Based on these analyses, we have chosen the RBF-FD approach, as it represents the best trade-off between accuracy and computational complexity. Furthermore this method is inherently present in the solution procedure, which is advantageous in two respects: (i) no additional approximation algorithms need to be implemented and (ii) the approximation error of the interpolated concentration is in the worst case of the same order of magnitude as the approximation of the partial differential operators involved~\cite{le2023guidelines}.

          %   Therefore, Shepard's interpolation method is excluded from further analysis due to its low accuracy. The accuracy of the stable RBF-FD and a more complex PU interpolation is less sensitive to the neighbourhood size $n$ and comparable for all sizes $n$ used. Based on these observations, we decide to use the RBF-FD interpolation method to obtain the solute concentration $u_\delta$ because it has good convergence behaviour, it is stable with respect to the template size and it has a with an accuracy comparable to that of the more complex PU approach and is favourable from the point of view of computational complexity.
          \begin{figure}[h]
              \centering
              \includegraphics[width=\textwidth]{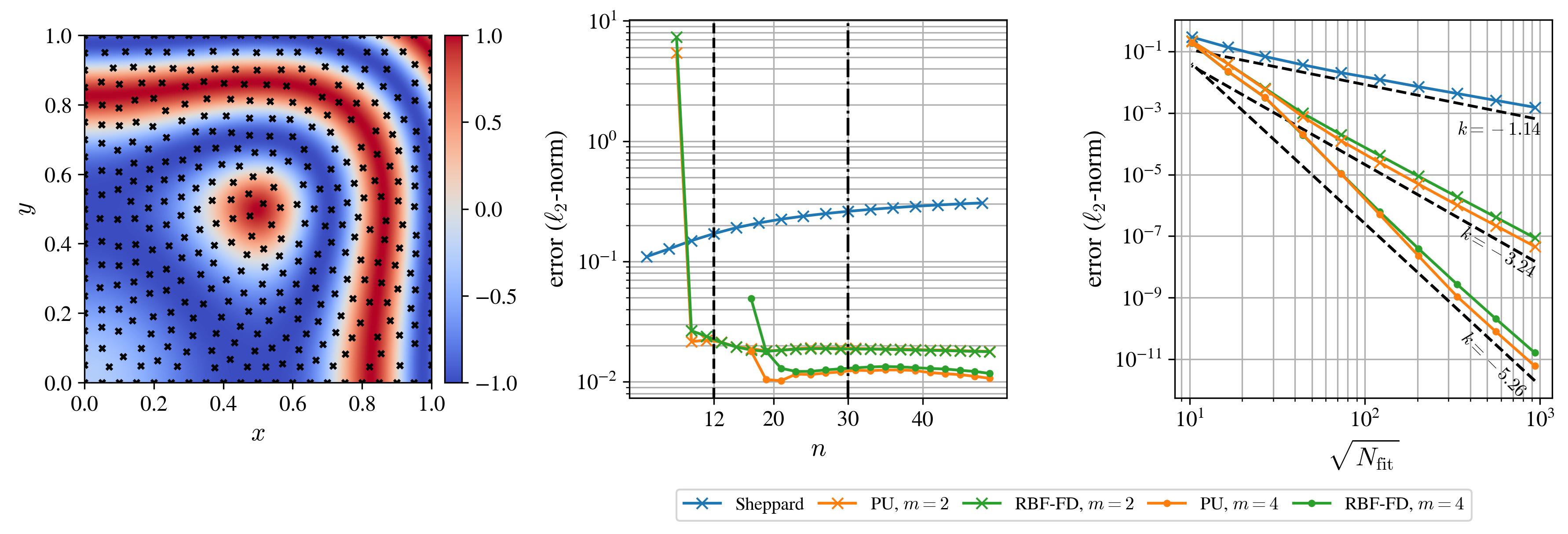}
              \caption{Analysis of interpolation performance on a two-dimensional Sibson's function. The left figure illustrates the interpolated field together with the fitting points marked as black crosses. The center figure shows the interpolation error with respect to the stencil size, $n$, in terms of the relative $\ell_2$ norm. The stencil sizes $n=12$ and $n=30$, recommended for $m=2$ and $m=4$, respectively, are shown with dashed vertical lines. The right figure shows the error with respect to the number of fitting points using recommended stencil sizes and $n=3$ for Sheppard's interpolation. The black dashed lines and their slopes, $k$, indicate the observed order of convergence.}
              \label{fig:interpolation}
          \end{figure}
\end{itemize}

\subsubsection{Grain envelope reconstruction}
\label{sec:reconstruction}
%After calculating the speed $v$ for each node $\x\ienv$ on the envelope, the nodes are translated to new positions following Equation~\eqref{eq:displacement} (see also line~\ref{alg:new_envelope} of Algorithm~\ref{alg:simulate}). However, the 
The spatial propagation of
%a set of 
the envelope nodes, as the grain grows,
%by a set of displacement vectors 
introduces two difficulties. First,
%if the envelope nodes are simply moved along the vectors perpendicular to the envelope boundary, 
the node spacing increases and the local internodal distance $h$ sooner or later violates the requirement of a quasi-uniform internodal distance for stable meshless approximations~\cite{wendland2004scattered,bayona2010rbf, slak2019generation}. Second, some of the nodes from the interior of the computational domain are eventually engulfed by the growing envelope and thus fall outside the domain. Both problems lead to an unstable solution procedure and must therefore be handled accordingly.

To avoid such behaviour, the entire domain is re-discretized when necessary (see line~\ref{alg:disp} of Algorithm~\ref{alg:simulate}), according to the criterion introduced in Section \ref{sec:discretization}. The first step of the re-discretization is a reconstruction of the envelope (see line~\ref{alg:re-discretize} of Algorithm~\ref{alg:simulate}) using a generalised surface reconstruction algorithm~\cite{surface}. This constructs a parametric representation of the envelope using splines $\gamma$ of the $k$-th degree over the translated nodes. Once the parametric representation of the envelope is obtained, a special surface DIVG algorithm~\cite{duh_surface} is used to populate it with a new set of nodes that are suitable for the RBF-FD approximation method~\cite{slak2019generation}.

The rest of the domain can then be discretized with the general DIVG node positioning algorithm.
%It is worth mentioning, that domain re-discretization is among the computationally most expensive steps in the solution procedure, but nevertheless required to assure accurate simulation of the dendritic growth.
%Readers interested in a detailed description of the surface reconstruction algorithm are to to its original proposal~\cite{surface}, nevertheless, a short description, including the Algorithm~\ref{alg:reconstruction}, is given below.

The surface reconstruction is performed by Algorithm~\ref{alg:reconstruction}, which takes a set of randomly ordered envelope nodes $\Gamma_\mathrm{e}$ about to be parametrized with a Jordan curve $\gamma: [a,b] \rightarrow \R^2$. To obtain the Jordan curve, the nodes are first ordered, therefore, a $k$-d tree with envelope nodes is constructed, which allows us to query for nearest neighbours, obtaining a permutation array, $\sigma$, in the process. To find the appropriate permutation $\sigma$, the list of ordered points is first initialized with an arbitrary starting point $\b x_0$. Its nearest neighbour $\b x_p$ is then obtained using the $k$-d tree and added to the list of ordered points. The remaining nodes are then inductively ordered: Once $\b x_j$ is in the list of ordered points, find its nearest neighbour $\b x _p$. If $\b x_p$ is not in the list of already ordered points, set $\b x_{j+1}$ to $\b x_p$, otherwise find the second (or $n$-th) nearest neighbour that is not yet in the list of ordered points. The process is repeated until all construction points from the envelope are ordered.

\begin{algorithm}
    \caption{Discretized boundary surface reconstruction.}
    \label{alg:reconstruction}
    \vspace{1mm}
    \textbf{Input:} Envelope boundary nodes $\Gamma_{e}$, nodal density function $h:~\Omega\to \R$, parametric spline degree $k$.\\
    \textbf{Output:} New domain discretization $\Omega^\star$. \\

    \begin{algorithmic}[1]

        \Function{discretize\_envelope}{$\Gamma_\mathrm{e}, h, k$}
        \State $\sigma \gets \Call{order\_nodes}{\Gamma_\mathrm{e}}$
        \Comment{Returns permutation of nodes $\sigma$.}
        \State $\gamma \gets \Call{fit\_spline}{\Gamma_\mathrm{e}, \sigma, k}$
        \Comment{Fits spline to ordered points.}
        \State $\Gamma_\mathrm{e}^\star \gets \Call{discretize}{\gamma, h}$
        \Comment{Obtaines new envelope discretization.}

        \State \Return $\Gamma_\mathrm{e}^\star$
        \EndFunction\\

        \Function{re-discretize\_domain}{$\Gamma_\mathrm{e}, h, k$}\label{alg:rediscretize}
        \State $\Gamma_\mathrm{e}^\star \gets \Call{discretize\_envelope}{\Gamma_\mathrm{e}, h, k}$\label{alg:re-discretize-envelope}
        \Comment{Obtaines new envelope discretization.}
        \State $\Omega^\star \gets \Call{discretize}{\Gamma_\mathrm{e}^\star, h}$\label{alg:discretize_rest}
        \Comment{Obtaines new domain discretization.}

        \State \Return $\Omega^\star$
        \EndFunction
    \end{algorithmic}
\end{algorithm}

\begin{figure}
    \centering
    \includegraphics[width=0.8\textwidth]{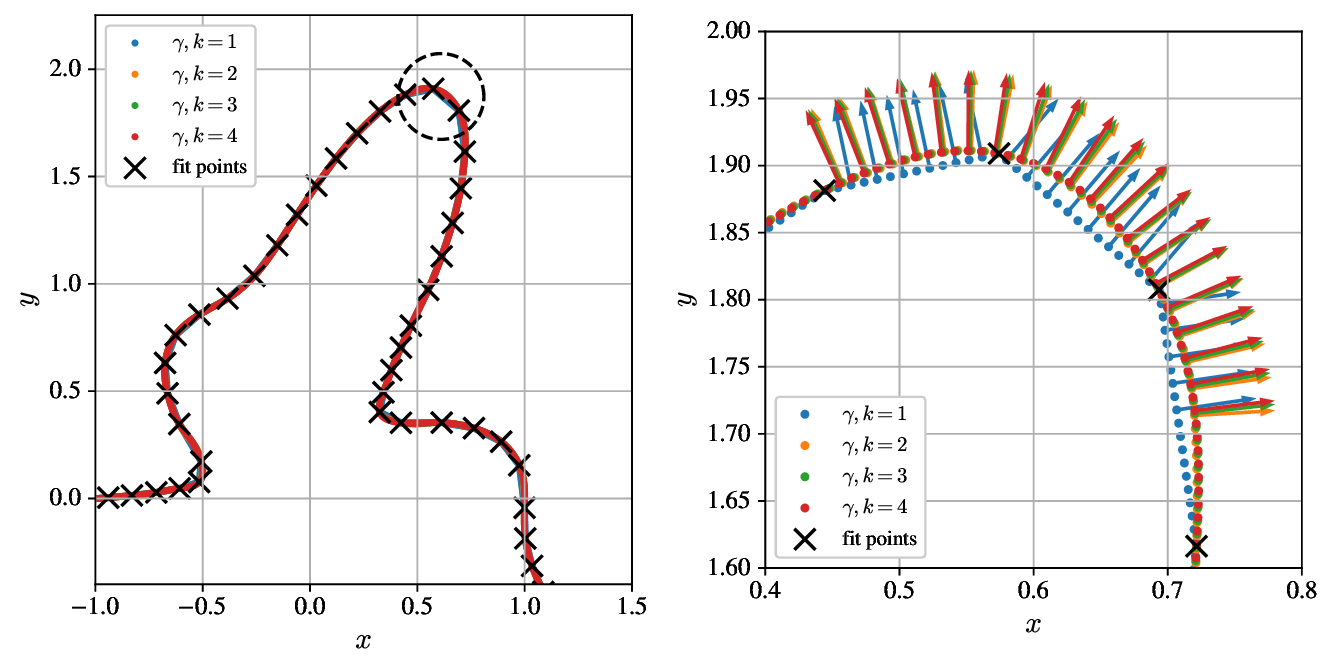}
    \caption{Example of surface reconstruction accompanied with normal vectors computed from the reconstructed spline $\gamma$ of degree $k\in \left \{1, 2, 3,4\right \}$. For clarity, the encircled area is zoomed in and shown on the right. We also show the normal directions in the right figure.}
    \label{fig:reconstruction}
\end{figure}

\hl{The ordered nodes are then used to re-construct the envelope shape by fitting a parametric spline} $\gamma$ \hl{of} $k$\hl{-th degree}. $\gamma$ \hl{is discretized with a given point spacing as the new discretized grain envelope,} $\Gamma_\mathrm{e}^\star$. \hl{Note that regular parametrization of} $\gamma$ \hl{is required to employ the surface DIVG}~\cite{duh_surface}. \hl{The surface DIVG takes a regular parametrization} $\gamma$, \hl{its Jacobian, a nodal spacing function} $h$ \hl{and a set of ``seed parameters'' from parameter space as input. The algorithm then returns a set of regularly distributed nodes on the curve, conforming to the spacing function} $h$. \hl{It is important to note that the node spacing does not take place directly in the target space} $\R^d$. \hl{Instead, the algorithm places parameters in the parametric space. The distance between two parameters} $\xi_1$ \hl{and} $\xi_2$ \hl{is then computed between points} $\gamma(\xi_1)$ \hl{and} $\gamma (\xi_2)$ \hl{and compared to} $h$, \hl{making the parameterization of} $\gamma$ \hl{irrelevant for assuring a good discretization quality as long as the determinant of the Jaccobian matrix is nonzero (see Equation~(2.6) in}~\cite{duh_surface}). For a graphical demonstration of the surface reconstruction algorithm, see Figure~\ref{fig:reconstruction}.

To discretize the rest of the domain (see line \ref{alg:discretize_rest} of Algorithm~\ref{alg:reconstruction}), an \texttt{inside check} algorithm is required to distinguish between the interior and exterior of $\Gamma_\mathrm{e}^\star$. To determine whether a node $\x$ is inside $\Gamma_\mathrm{e}^\star$, a scalar product
\begin{equation}
    \label{eq:insidecheckgamma}
    \left \langle \x - \gamma(t_{\min}), c \b n(t_{\min}) \right \rangle,
\end{equation}
is computed. Here, $\gamma(t_{\min})$ is a node on the envelope closest to $\x$ \hl{and} $\b n(t_{\min})$ \hl{is the normal vector in the same node. The normal vector is computed through a tangent vector} $\b t(t_{\min}) = (t_x, t_y)  \propto {\gamma}'(t_{\min})$ \hl{yielding} $ \b n(t_{\min}) = (-t_y, t_x)$. If the scalar product~\eqref{eq:insidecheckgamma} is negative, $\x$ is in the interior of $\Gamma_\mathrm{e}^\star$, otherwise not. The constant
\begin{equation}
    \label{eq:constant}
    c = -\mathrm{sgn}(\left \langle  \x_\mathrm{int} - \gamma(t_{\min}), \b n (t_{\min})\right \rangle)
\end{equation}
is computed for an arbitrary node $\x_\mathrm{int}$ from the interior and ensures that the normals point outwards. In Equation~\eqref{eq:constant} $t_{\min}$ is determined so that the $\gamma(t_{\min})$ is closest to $\x_\mathrm{int}$.

\emph{Note:} The reconstruction algorithm assumes the nodes from $\Gamma_\mathrm{e}$ must be \emph{dense enough} to adequately describe the curve~\cite{surface}. In our particular case of dendritic growth simulation we must
%take care of this problem by using $h$-refinement with increasing density of nodes towards the envelope, assuring the internodal distance is much smaller than the resolution of the characteristic length of envelope details.
choose a node spacing on and in the vicinity of the envelope, such that the internodal distance is much smaller than the radius of curvature of the envelope.

\subsubsection{Choice of the spline degree for accurate envelope reconstruction}
%\label{sec:spline}
The effect of the spline degree is demonstrated in Figure~\ref{fig:reconstruction} for degrees $k\in \left \{1, 2, 3,4\right \}$. As expected, a linear interpolation between two neighbouring nodes ($k=1$) does not yield a satisfying representation of the envelope. %Additionally, the normal vectors (shown in Figure~\ref{fig:reconstruction} on the right) computed from $\gamma(k=1)$ are also inadequate.

\begin{figure}[h]
    \centering
    \includegraphics[width=0.9\linewidth]{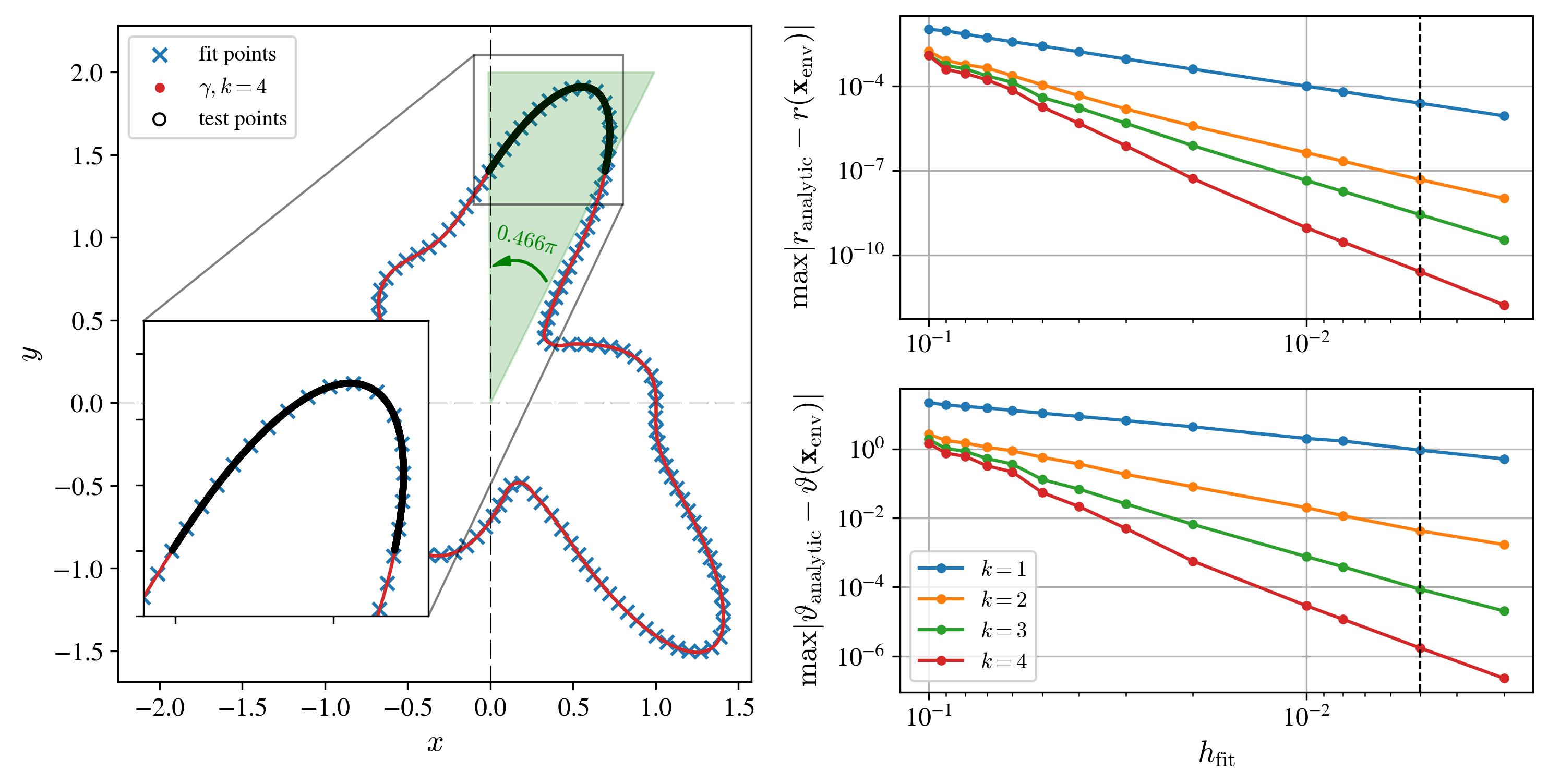}
    \caption{Parametric domain used to evaluate the accuracy of surface reconstruction (left) and reconstruction quality analysis (right) in terms of maximum radius error (top right) and maximum normal angle error (bottom right).}
    \label{fig:reconstruction_accuracy}
\end{figure}

Next, a more thorough discussion on the optimal spline degree $k$ for accurate surface reconstruction is given. The reconstruction quality is evaluated on a parametric function
\begin{align}
    \label{eq:parametric_test}
    R(t)    & =\lvert \cos(1.5t) \rvert ^{\sin(3t)} \\
    \b r(t) & =R(t)(\cos(t),\sin(t))
\end{align}
for values $t$ corresponding to the dendrite-tip-like shape (see Figure~\ref{fig:reconstruction_accuracy} left
for clarity).

A set of equally spaced fit nodes is taken to obtain a parametric representation of the boundary. The internodal distance $h_\mathrm{fit}$ of fit points ranges from $0.1$ to $0.003$, resulting in $124$ and $4094$ fit points respectively. The test points, where the reconstruction quality is evaluated, are generated by discretizing the parametric representation of the boundary $\gamma$ with a uniform internodal distance $h_\mathrm{test}=0.005$, yielding approximately $2450$ test nodes.

%Next, we evaluate the effect of the spline degree,  on the envelope reconstruction quality (line~\ref{alg:re-discretize-envelope} of Algorithm~\ref{alg:reconstruction}). 

The quality of the surface reconstruction is evaluated in terms of (i) the maximum error of the point position in terms of the distance $r=\sqrt{x^2+y^2}$ (top right in Figure~\ref{fig:reconstruction_accuracy}) and (ii) the maximum error of the angle $\vartheta$ of the normal vector $\b n$ (bottom right in Figure~\ref{fig:reconstruction_accuracy}); both within the domain of interest marked in light green in Figure~\ref{fig:reconstruction_accuracy} (left). From the analysis we conclude that the parametric representation of the envelope can be used to obtain an accurate reconstruction of the dendritic grain envelope, and that the spline degree $k$ significantly affects the accuracy of the reconstruction.

Note that a dashed vertical line has been added to the top and bottom right plots in Figure~\ref{fig:reconstruction_accuracy} to mark the internodal distance at which the density of the fit points is approximately equal to the density of the reconstructed nodes, i.e.\ $h_\mathrm{fit} \approx h_\mathrm{test}$. This is also mostly the case in the dendritic growth simulation governed by the GEM.

%In terms of accuracy, we confirm our previous observations that reconstruction with spline degree $k=1$ yields the least accurate results in terms of both distance to the origin $r(t)$ and normal vector angle $\vartheta(t)$. 

From the presented results, it is not possible to objectively deduce which spline degree $k$ yields a sufficiently accurate parametric representation of the dendritic grain envelope. The safest option is to choose the highest proposed spline degree, i.e.\ $k=4$, especially considering that the additional computational cost due to higher spline degrees is marginal.

A more systematic discussion is continued in Section~\ref{sec:istropy}, where the final decision on the spline degree $k=2$ used throughout the rest of this work is made.

\subsection{Implementation remarks}
The entire MIT solution procedure from Algorithm~\ref{alg:simulate} using meshless methods is implemented in C++ environment. The projects implementation\footnote{Source code is available at \url{https://gitlab.com/e62Lab/public/2022_p_2d_dendritic_growth} under tag \emph{v1.1}.\label{git}} is strongly dependent on our in-house development of a meshless C++ \emph{Medusa library}~\cite{slak2021medusa}.

The code was compiled using \texttt{g++ (GCC) 11.3.0 for Linux} with \texttt{-O3 -DNDEBUG} flags on \texttt{AMD EPYC 7702 64-Core Processor} computer. OpenMP API has been used to run parts of the solution procedure in parallel on shared memory. %due to execution timings, the CPU frequency was fixed to 2.27 GHz with disabled boost functionality and assured CPU affinity using the \texttt{taskset} command. 
Post-processing was done using Python 3.10.6 and Jupyter notebooks, also available in the provided git repository\footref{git}.

\section{MIT solution procedure verification: Isotropic growth of grain envelope}
\label{sec:istropy}

The goal of this section is analyse how the behaviour of the full MIT procedure depends (i) on the order of the RBF-FD approximation and (ii) on the order of the spline reconstruction of the envelope, and to quantify the errors linked to these methods.

\begin{figure}[h]
    \centering
    \begin{subfigure}[b]{0.425\textwidth}
        \centering
        \includegraphics[width=\linewidth]{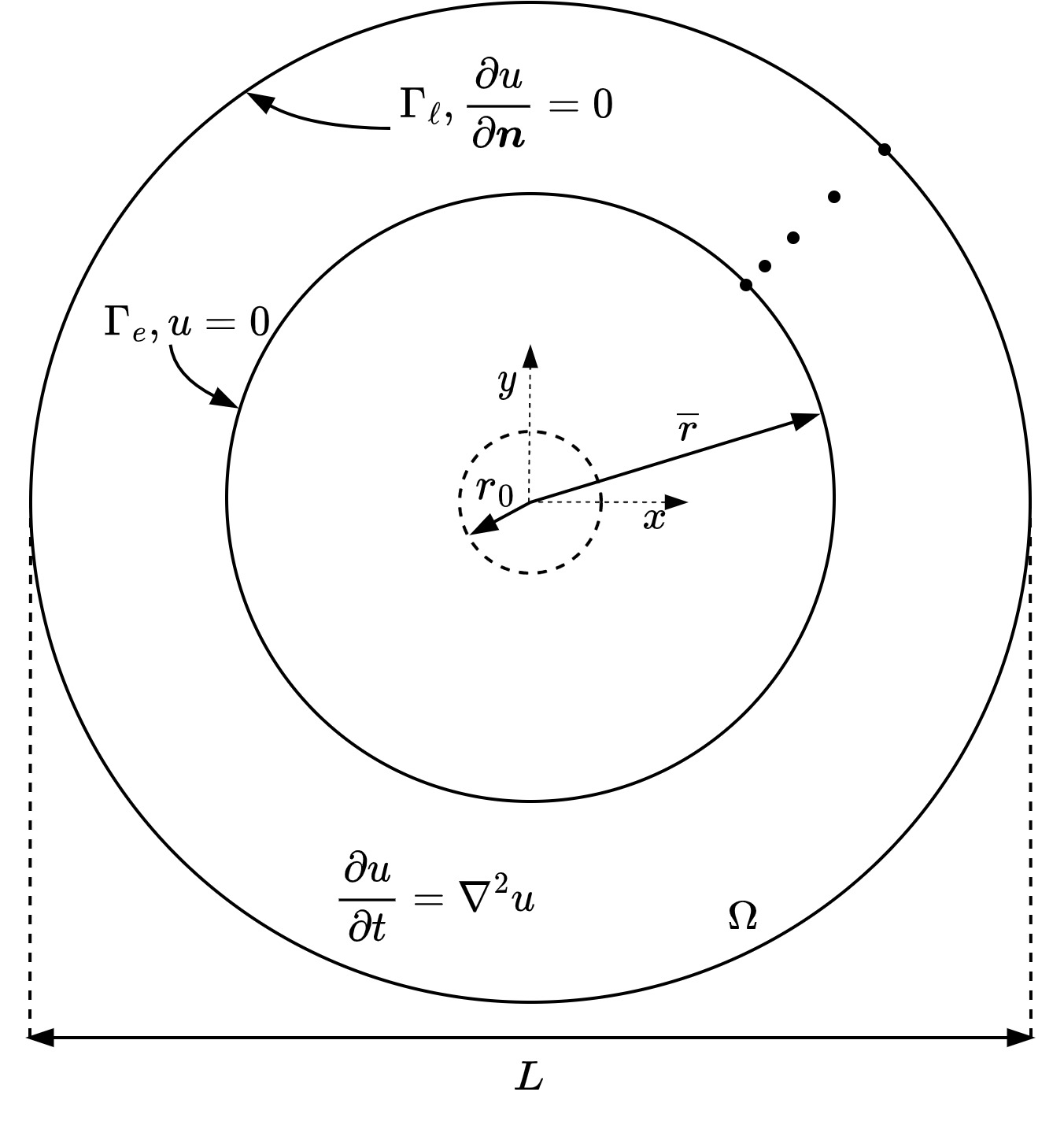}
        \caption{Schematic representation of the isotropic problem. In the top right corner a set of black nodes is used to illustrate the \h-refinement towards the envelope.}
        \label{fig:sketch_isotropic}
    \end{subfigure}
    \hfill
    \begin{subfigure}[b]{0.49\textwidth}
        \centering
        \includegraphics[width=\linewidth]{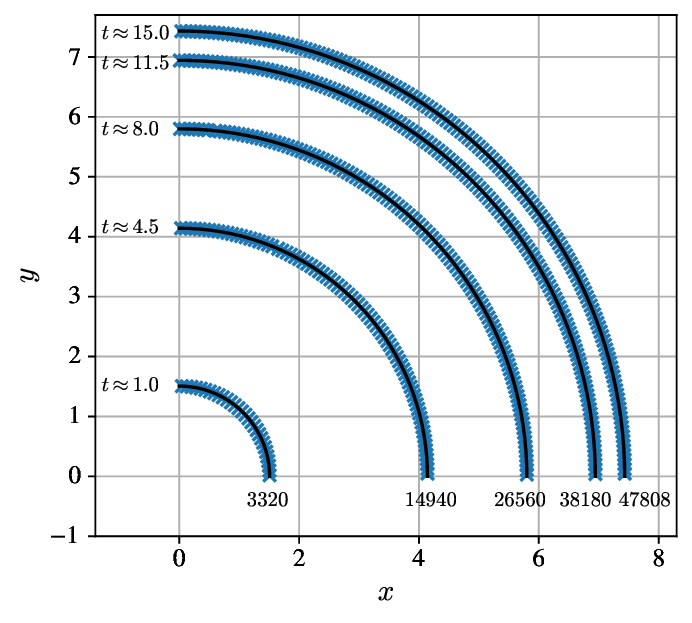}
        \caption{Isotropic growth time lapse for 5 simulations times $t$ with the total number of time steps written at bottom of each envelope.}
        \label{fig:isotropy}
    \end{subfigure}
    \caption{Presentation of the isotropic growth case.}
    \label{fig:isotropy_presentation}
\end{figure}

To verify the full MIT solution procedure we propose a simple isotropic growth test, effectively forcing the Equation~\eqref{eq:displacement} to take the following form
\begin{equation}
    \x_i' = \x_i + v(u_{\delta,i}) \b n_i \mathrm{d}t, \quad \forall \x_i \in \Gamma_\mathrm{e},
\end{equation}
and thus discarding the preferred directions of growth. With such setup, the grain envelope is expected to maintain its initial circular shape through all computational time steps and can therefore be used for verification of the method. The test is shown schematically in Figure~\ref{fig:sketch_isotropic}. The computational domain is a circle of diameter $L=20$ with Neumann boundary conditions enforced on the outer (liquid) boundary $\Gamma_\ell$ and a Dirichlet boundary condition of $u=0$ enforced on the growing grain envelope $\Gamma_\mathrm{e}$. At time $t=0$, the envelope is initialized with a circle of radius $r_0 = 0.22239$ and the liquid in the domain $\Omega$ has a uniform initial concentration of $u_0=0.18$.

\begin{wrapfigure}{r}{0.5\textwidth}
    \vspace{-0.77cm}
    \centering
    \includegraphics[width=0.49\textwidth]{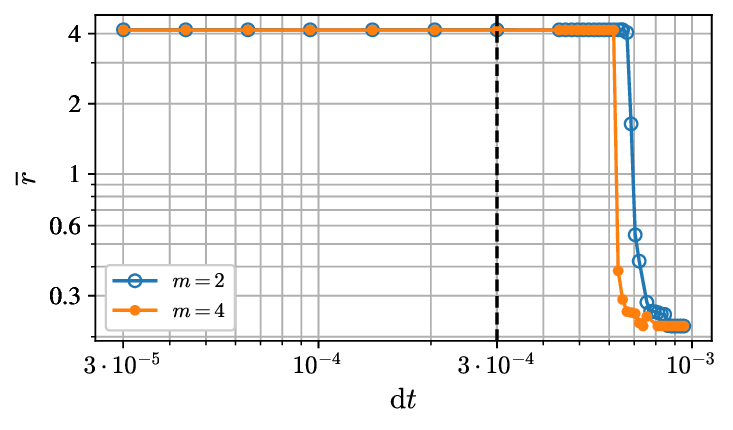}
    \caption{Average radius of the isotropic envelope at simulation time $t=4.5$ with respect to time step $\mathrm{d}t$.}
    \label{fig:timestep_conv}
\end{wrapfigure}
In Figure~\ref{fig:isotropy}, a time-lapse of the isotropic envelope calculated with MIT shows the circular envelope expansion. For clarity, a black line is constructed to demonstrate an \emph{exact} circular shape, with the radius equal to the mean distance of the envelope nodes from the grain center.

Before we perform the analysis of MIT spatial discretization error, a scan over the time step size is performed to ensure an appropriate temporal discretization. \hl{The stability of the method in }$h$\hl{-refined node arrangement is limited to the worst case, i.e. to the smallest internodal distance, therefore the optimal approach would be to scale time-step according to the varying }$h$\hl{. In a recent paper}~\cite{Dobravec2023} \hl{authors used quad-tree time step adaptivity to take this into account. In this paper, however, we use a global time step, which we set according to the stability criterion for the smallest internodal distance with additional safety margin. To ensure we indeed use a small enough time step, we also performed numerical analysis presented in} Figure~\ref{fig:timestep_conv}, where the average envelope radius $\overline{r}$ at simulation time $t=4.5$ is shown with respect to the time step $\mathrm{d}t$.
%According to the Courant stability criteria,
Too large time steps ($\mathrm{d}t\ge 6\cdot 10^{-4}$) result in divergent solutions, while time steps $\mathrm{d}t\le 3\cdot 10^{-4}$ give stable solutions for both the second and fourth order discretization approaches. Therefore, we use $dt = 3\cdot 10^{-4}$ in all following discussions.

\begin{figure}[h]
    \centering
    \includegraphics[width=0.85\textwidth]{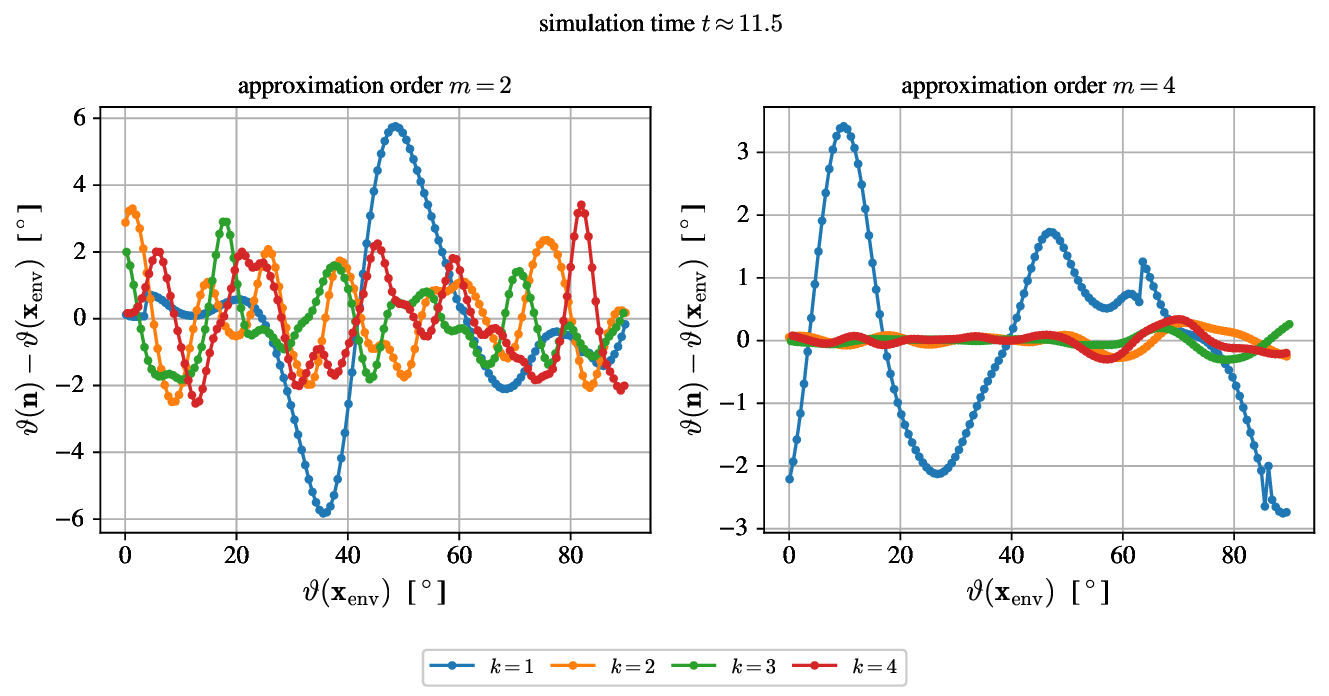}
    \caption{Shape distortion of the isotropic envelope: the difference in angle between the computed normal $\vartheta (\b n)$ and the angle of the radial vector pointed towards the envelope position  $\vartheta (\b x_{\rm env})$ at $t\approx 11.5$ using second order RBF-FD (left) and fourth order RBF-FD (right). $k$ is the order of the spline used for envelope reconstruction.}
    \label{fig:isotropy_angles}
\end{figure}

%First, we will observe discrepancy between expected circular shape and shape we get after different times of MIT simulation. 
Next, Figure~\ref{fig:isotropy_angles} shows the
%difference between the angle of calculated normal $\vartheta (\b n)$ and the angle of the radial vector pointed towards the envelope position $\vartheta (\b x_{\rm env})$ at simulation time $t\approx 11.5$. 
distortion of the calculated normal angle, $\vartheta (\b n)$, as function of the radial angle, $\vartheta (\b x_\mathrm{env})$.
The deviations from the ideal circular shape do not significantly depend on the spline degree, $k$, for $k > 1$, while $k=1$ is too low for adequate surface reconstruction. We also observe a significant improvement in MIT accuracy when using a higher order RBF-FD (compare the left and right plots of Figure~\ref{fig:isotropy_angles}, corresponding to the second and fourth order RBF-FD approximation, respectively). The shape deviations are presented in a more quantitative manner in Table~\ref{tab:tabular}, where the extrema of the distortion in the envelope radius, from the mean radius, $\overline{r}$, and the largest distortions of the normal angles, $\vartheta (\b n) - \vartheta (\b x_\mathrm{env})$, are tabulated.

\begin{table}[H]
    \begin{center}
        \begin{tabular}{  l |c c c c c c}
            \toprule
                                             &                 & $\overline{r}$ & $\overline{r} - \min(r)$ & $\max(r) - \overline{r} $ & $\max(\vartheta(\b n)  -\vartheta(\x_\mathrm{env}))$ [$^\circ$] & $\min(\vartheta(\b n)  -\vartheta(\x_\mathrm{env}))$ [$^\circ$] \\ \midrule

            \multirow[0]{4}{*}{\rot{$m= 2$}} & $\gamma(k = 1)$ & $7.3690$       & $0.0340$                 & $0.0357$                  & $2.9411$                                                        & $-3.0012$                                                       \\
                                             & $\gamma(k = 2)$ & $7.3967$       & $0.0285$                 & $0.0203$                  & $3.9000$                                                        & $-4.0700$                                                       \\
                                             & $\gamma(k = 3)$ & $7.3970$       & $0.0225$                 & $0.0301$                  & $2.9738$                                                        & $-4.4241$                                                       \\
                                             & $\gamma(k = 4)$ & $7.3983$       & $0.0263$                 & $0.0238$                  & $2.8899$                                                        & $-3.0776$                                                       \\\midrule
                                             & $\gamma(k = 1)$ & $7.3514$       & $0.0322$                 & $0.0321$                  & $3.2313$                                                        & $-2.8206$                                                       \\
                                             & $\gamma(k = 2)$ & $7.3913$       & $0.0037$                 & $0.0031$                  & $0.3517$                                                        & $-0.2851$                                                       \\
                                             & $\gamma(k = 3)$ & $7.3908$       & $0.0036$                 & $0.0021$                  & $0.2039$                                                        & $-0.2423$                                                       \\
            \multirow[0]{-4}{*}{\rot{$m=4$}} & $\gamma(k = 4)$ & $7.3914$       & $0.0014$                 & $0.0032$                  & $0.2227$                                                        & $-0.3810$                                                       \\
            \bottomrule
        \end{tabular}
    \end{center}
    \caption{Shape distortion from a circular shape at $t=11.5$ depending on the order of the RBF-FD approximation, $m$, and the on order of the splines for surface reconstruction, $k$. The distortion is represented in terms of the largest deviations from the mean radius, $r - \overline{r}$, and the largest distortions of the normal angles, $\vartheta (\b n) - \vartheta (\b x_\mathrm{env})$.}
    \label{tab:tabular}
\end{table}

From the above analysis we can conclude that: (i) splines of order $k > 1$ should be used for envelope reconstruction to achieve good accuracy, (ii) shape distortions depend on the order of the RBF-FD and high-order ($m=4$) is required for high accuracy of the envelope normal direction, while the accuracy of the mean size (radius) is satisfactory even at low RBF-FD orders ($m=2$).
%accurate calculation of the in terms of mean radius $\overline{r}$ there is no obvious advantage in using a high order RBF-FD, especially when considering the computational cost, however, an advantage can be seen in smaller extremal values of normal vector angle errors.

\begin{figure}[h]
    \centering
    \includegraphics[width=\linewidth]{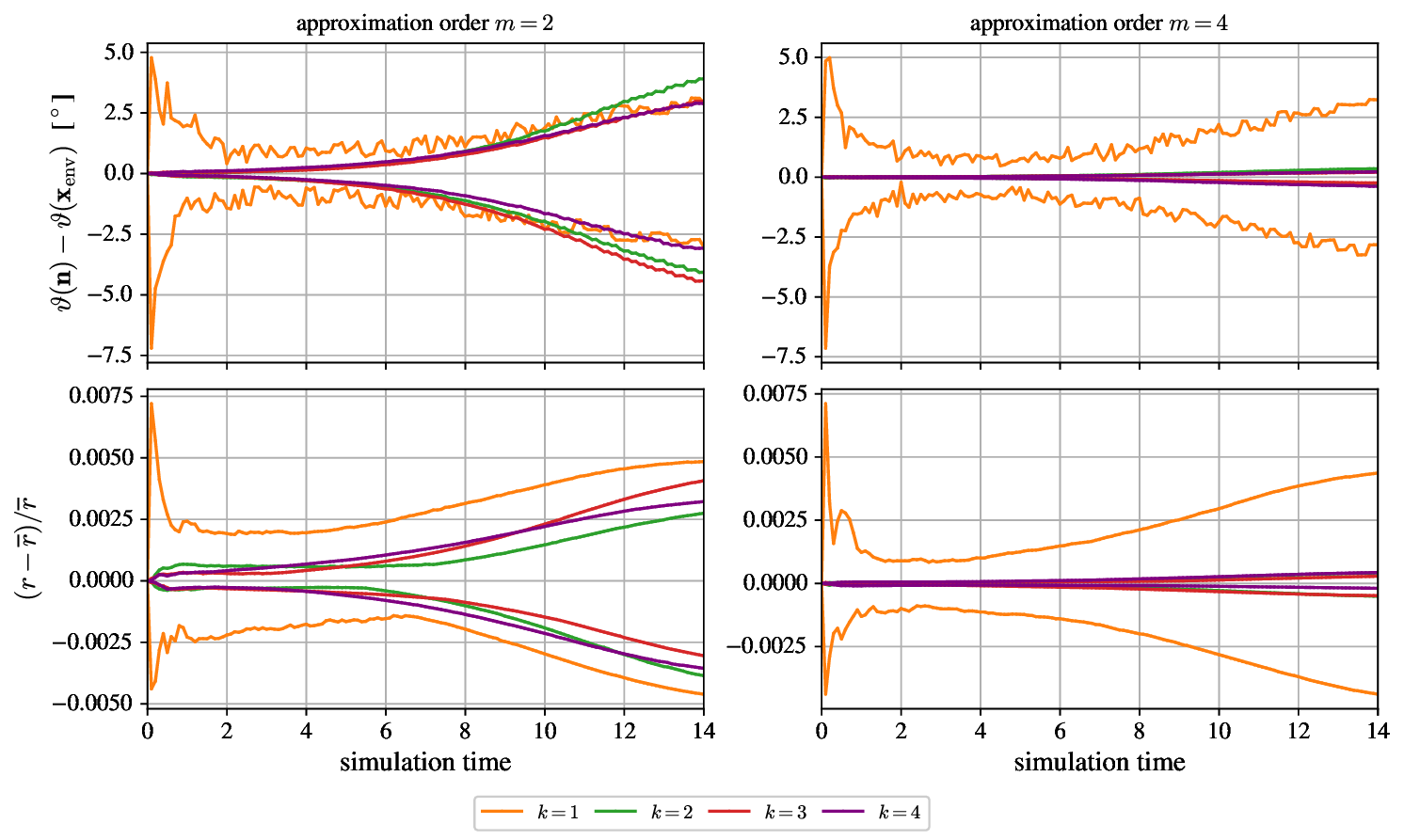}
    \caption{Time evolution of the shape distortion of the isotropic envelope. The top row shows distortion in normal vector angles and the bottom row shows distortion in terms of radius, the left column shows second order and the right column shows fourth order RBF-FD approximation.}
    \label{fig:isotropy_timelapse}
\end{figure}

The temporal evolution of the two error metrics ($\vartheta (\b n)$ - $\vartheta (\b x_{\rm env})$ and $\overline{r} \pm [\mathrm{min}(r), \mathrm{max}(r)]$) is assessed in Figure~\ref{fig:isotropy_timelapse}. Except for $k=1$ splines, the error strictly increases with time. It is now even more evident that $k=1$ does not give a sufficiently accurate surface reconstruction. The differences between the second- and fourth-order RBF-FD approximations are relatively small (less than 0.5\% difference in the radius even at the end of the simulation), thus the use of the computationally much more demanding fourth-order method might not be justified.

\section{Simulation of dendritic growth}
\label{sec:results}
%The simulation is started by building the initial domain discretization, where the initial domain is a difference between the square domain with dimensionless side length $L = 20$ representing the outer liquid boundary and the initial grain envelope represented as a circle with dimensionless radius $r_0 = 0.22239$. The node density is prescribed to linearly decrease from $\he$ on the dendrite envelope towards the $h_\ell = 0.1$ on the outer liquid boundary. We use mixed boundary conditions with Neumann boundary condition on the outer liquid boundary ($\Gamma_\ell$) and Dirichlet boundary condition on the dendritic envelope ($\Gamma_\mathrm{e}$) with prescribed solute concentration. The initial solute concentration is homogeneous and equals $\Omega_0 = 0.18$, except on the dendritic envelope where the concentration is constant throughout the simulation and equals 0. A graphical sketch of the problem including boundary conditions is given in Figure~\ref{fig:sketch}. For clarity, the initial dendritic envelope is sketched with a dashed circle. Additionally, a sketch of the linearly \h-refined discretization nodes from the outer liquid boundary towards the envelope is demonstrated with black nodes in the top right corner of the sketch. Note that all our simulations are stopped after the maximum tip position reaches $x_{\text{tip}}=8.9$. Of course, other stopping criteria could be used, but this was sufficient for our analyses.

\begin{wrapfigure}{r}{0.4\textwidth}
	% \vspace{-0.5cm}
	\vspace{-2.5cm}
	\begin{center}
		\includegraphics[width=0.99\linewidth]{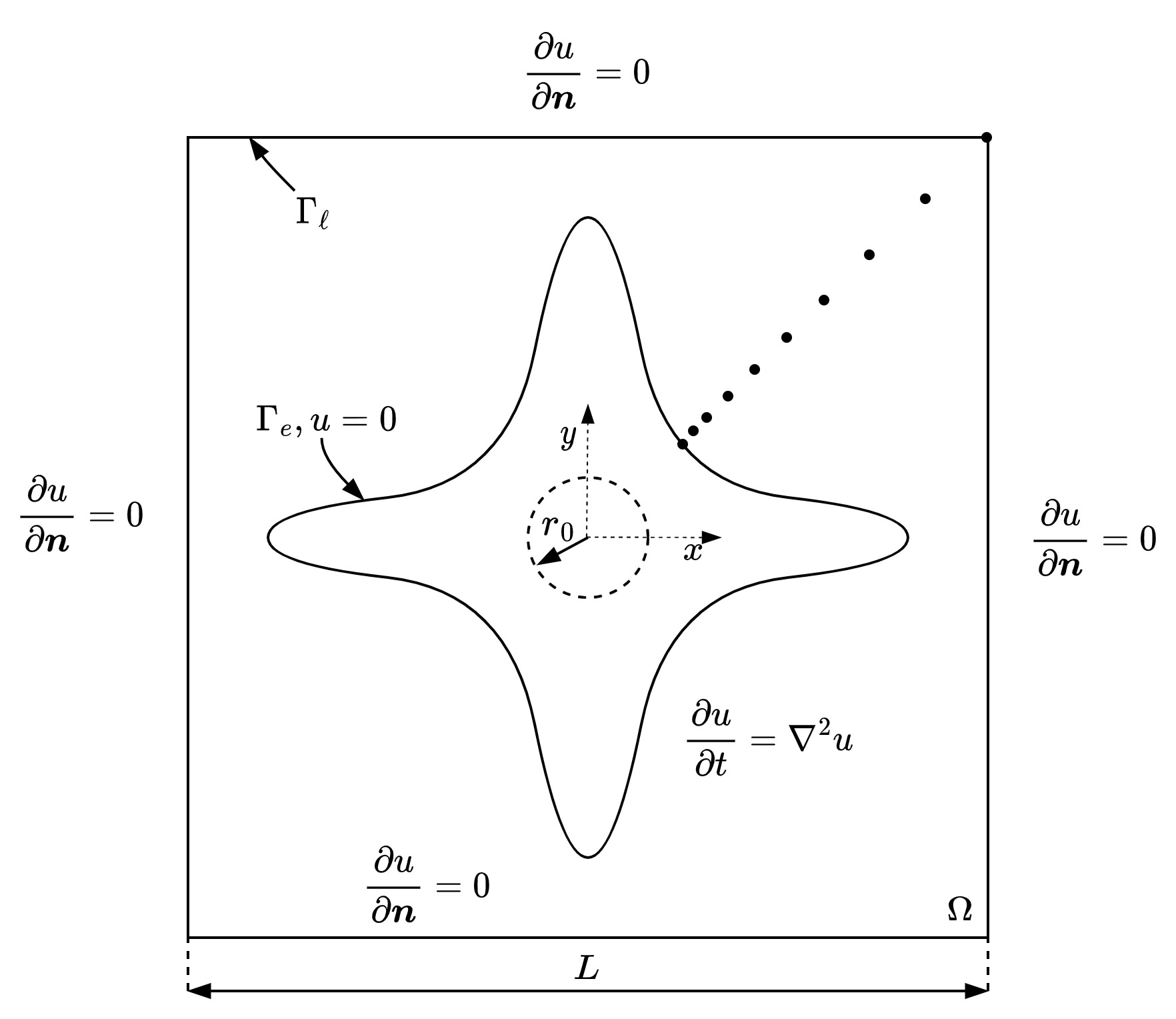}
		\caption{Schematic representation of dendritic growth problem.}
		\label{fig:sketch}
	\end{center}
\end{wrapfigure}
%In this section, dendritic growth simulations employing the standard PFIC-based and the proposed MIT-based GEM are studied. Their performance in terms of accuracy is analysed.
In this section we apply the MIT to the full GEM model to perform simulations of a dendritic grain and we compare the results to the PFIC method. The dendrite growth problem is schematically shown in Figure~\ref{fig:sketch}. The computational domain is a square of side length $L=20$ with Neumann boundary condition enforced on the outer (liquid) boundary $\Gamma_\ell$ and a Dirichlet boundary condition of $u=0$ enforced on the growing dendrite envelope $\Gamma_\mathrm{e}$. At time $t=0$, the dendrite envelope is initialized with a circle of radius $r_0 = 0.22239$ and the liquid in the domain $\Omega$ has a uniform initial concentration of $u_0=0.18$. The simulations are ran until the primary dendrite tip growing along the $x$ axis reaches position $x_{\text{tip}}=8.9$.
%The time step is a function of node spacing on the envelope $\he$, i.e., $\mathrm{d}t = 0.1 \frac{\he^2}{2}$.

The node spacing on the envelope, $h_\mathrm{e}$ is kept constant throughout the simulation as described in Section~\ref{sec:reconstruction}. In the domain, the inner nodes are generated such that the node spacing increases linearly from $h_\mathrm{e}$ on the dendrite envelope to $h_\ell = 0.1$ on the outer liquid boundary, according to the nodal density function defined in Equation~\eqref{eq:slakNodes}. Note that $h_\ell = 0.1$ is constant and identical for all simulations, regardless of $\he$. We use the RBF-FD approximation method, where the approximation basis consists of cubic PHS $r^3$ and monomials of second, $m=2$, and fourth, $m=4$, degree. In all cases, the differential operator approximation and the local interpolation of the concentration to the stagnant film are performed with the same RBF-FD setup. We present the results for both, i.e., the second and fourth order RBF approximations.

\begin{figure}[H]
	\centering
	\includegraphics[width=0.75\textwidth]{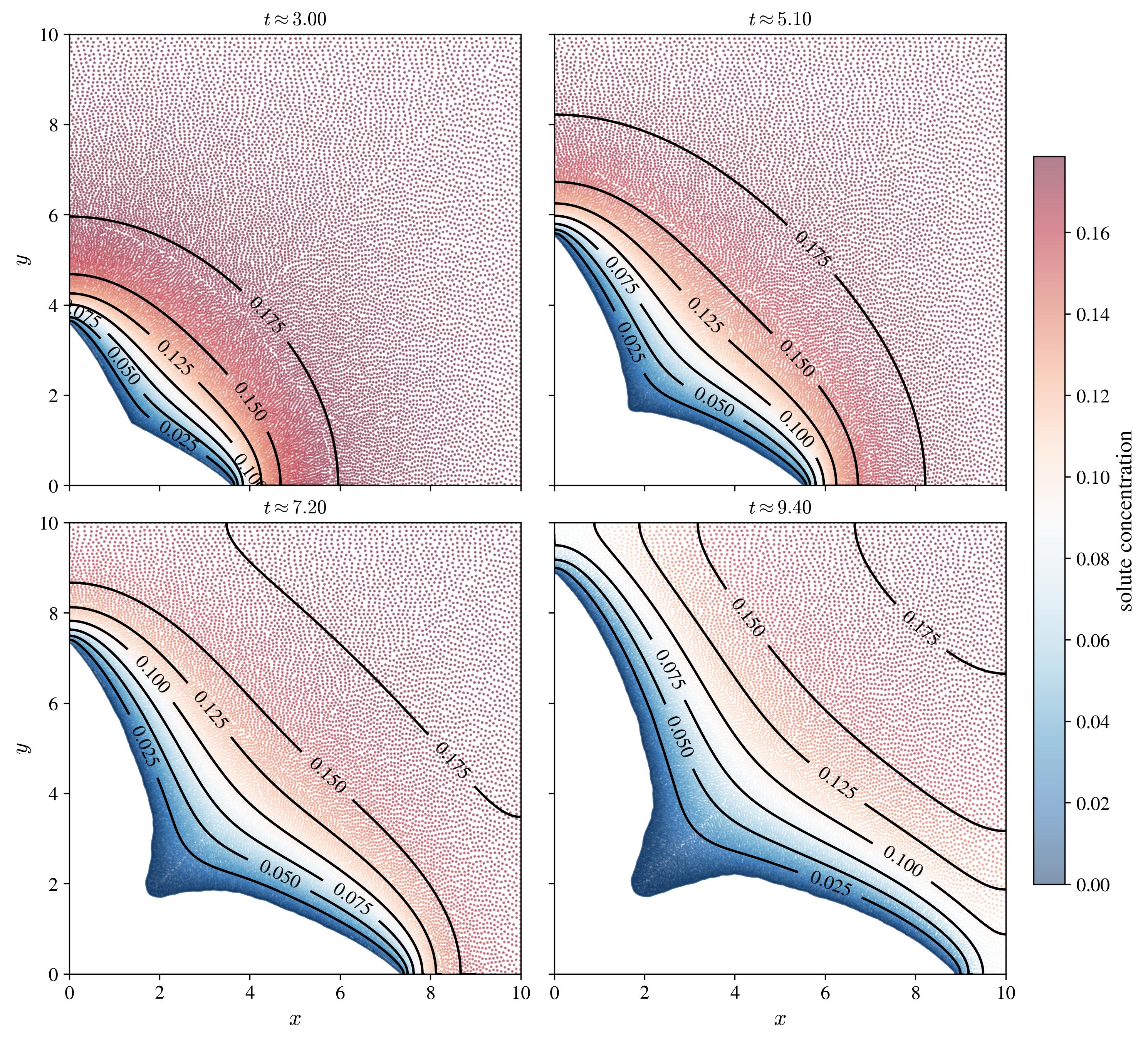}
	\caption{Growth of the dendritic grain envelope, the concentration field in the liquid and the adaptive discretization for $\he=0.028$}
	\label{fig:solution_fields}
\end{figure}

Figure~\ref{fig:solution_fields} shows the solute concentration fields for selected simulation time steps computed with MIT using the internodal spacing $\he \approx 0.028$ on the dendrite envelope, effectively resulting in $67\,544$ and $77\,168$ discretization nodes at the first and last time steps, respectively. Note that due to the symmetry of the problem, only the first quadrant of the domain $\Omega$ is shown, although the entire dendrite was simulated.

\begin{figure}[h]
	\centering
	\centering
	\begin{subfigure}[b]{0.47\textwidth}
		\centering
		\includegraphics[width=\textwidth]{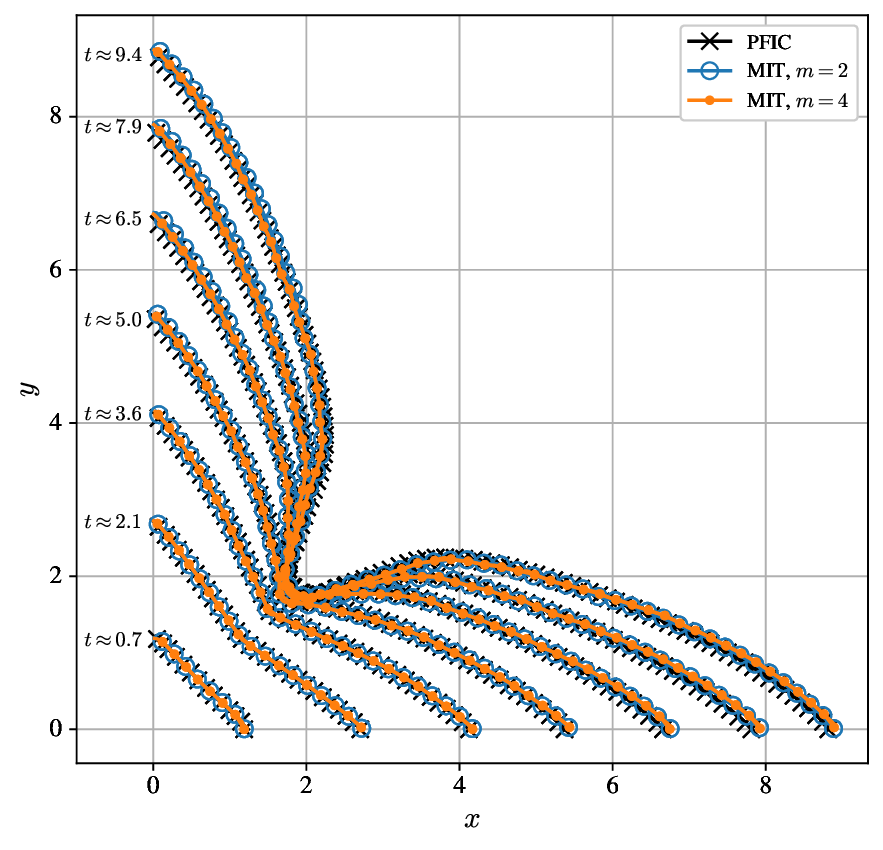}
		\caption{Dendritic envelope shape time lapse.}
		\label{fig:envelop_shapes}
	\end{subfigure}
	\hfill
	\begin{subfigure}[b]{0.515\textwidth}
		\centering
		\includegraphics[width=\textwidth]{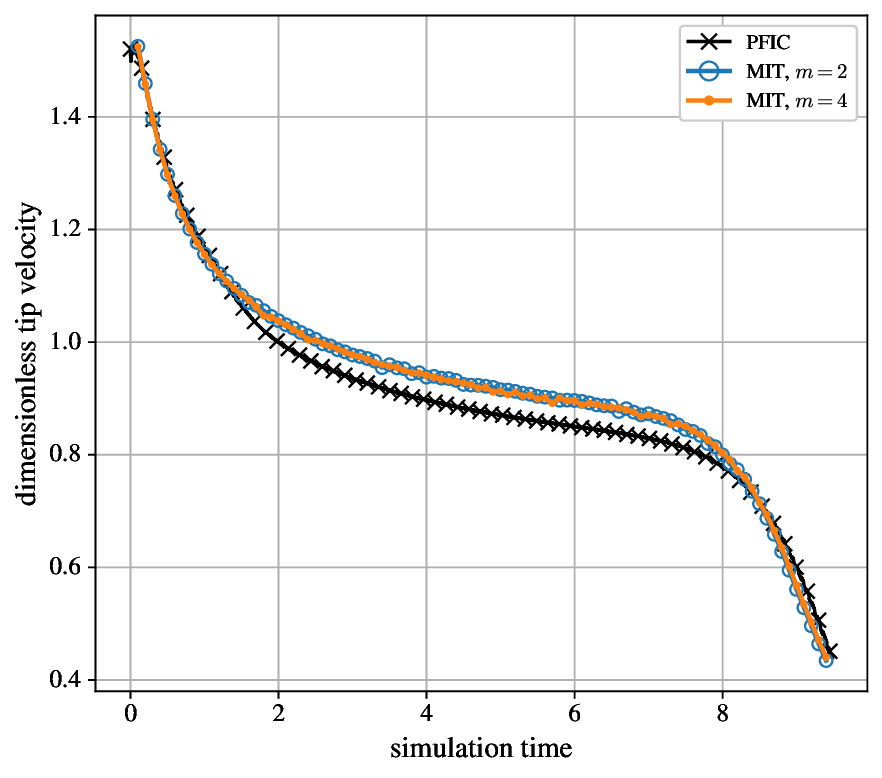}
		\caption{Primary dendrite tip velocity with respect to simulation time.}
		\label{fig:conv}
	\end{subfigure}
	\caption{Comparison of MIT and PFIC on dendritic growth problem for $\he\approx 0.028$.}
	\label{fig:dendritic_envelopes}
\end{figure}

To compare the MIT and PFIC methods, the envelope ($\Gamma_\mathrm{e}$) is compared at different times in Figure~\ref{fig:envelop_shapes}, along with the primary dendrite tip velocity with respect to time $t$ in Figure~\ref{fig:conv}. Although the results of both methods are generally consistent and are very close, there are some differences that are the result of the differences of the two solution approaches. Nevertheless, the first remark can be made: Since PFIC and MIT generally agree, MIT can also be used to resolve the GEM.

%Again, the MIT and PFIC in general agree well-based approach results in numerical solutions that are in good agreement with the PFIC-based approach. The computed dimensionless tip velocity matches well with the reference PFIC-based solution within the entire range of tip positions. Furthermore, convergent behaviour with respect to discretization quality of the dendrite envelope and its neighbourhood can also be observed. The green color is used to denote the coarsest discretization, while the blue color is used to denote the finest discretization.

\subsection{Anisotropy induced by the scattered node arrangement}
\begin{figure}[H]
	\centering
	\includegraphics[width=\textwidth]{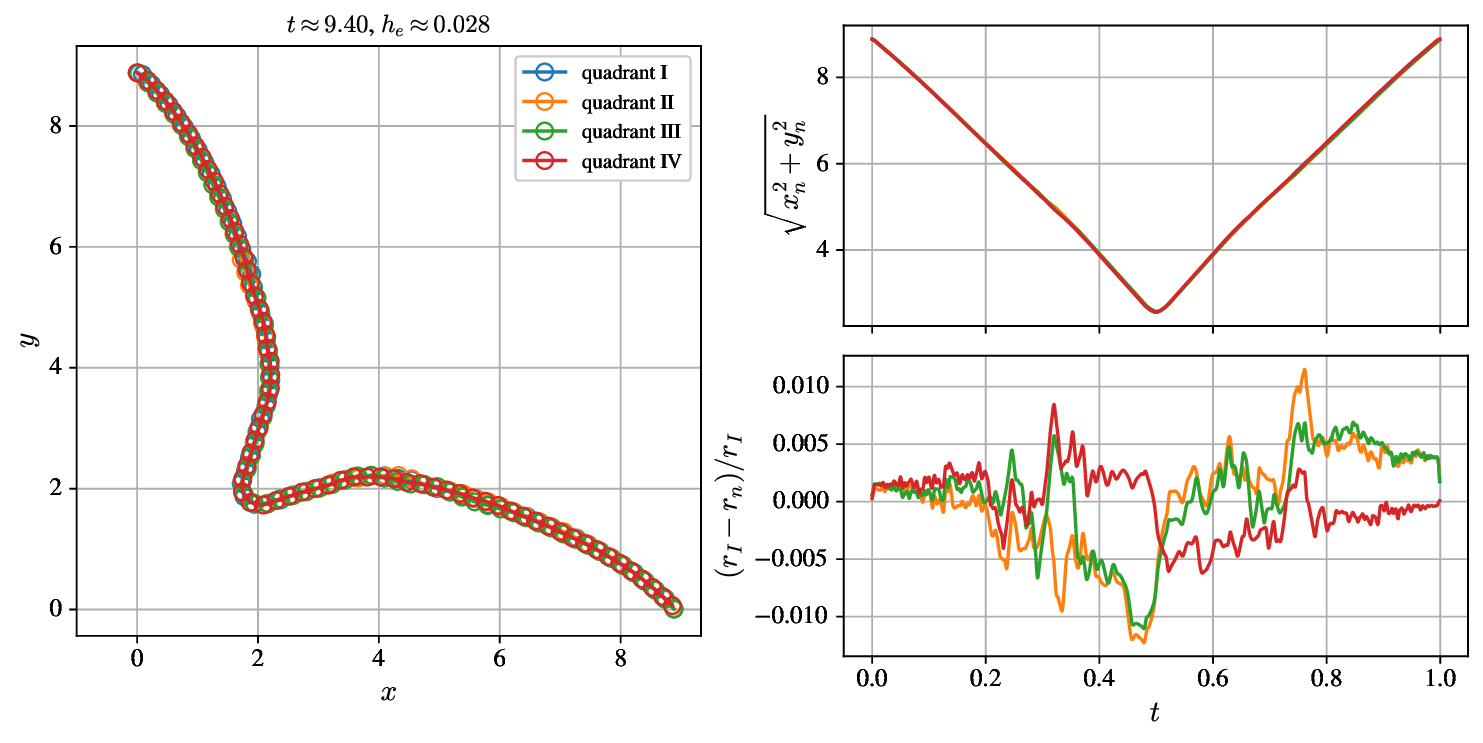}
	\caption{Analysis of anisotropy induced by the node arrangement. On the left, we show the envelope parts mapped to the first quadrant, on the right the evaluation process is presented.}
	\label{fig:induced_anisotropy}
\end{figure}

Although the growth of the dendritic envelope is symmetric according to the model, we can expect irregularities in this respect due to the numerical errors. This is even more pronounced in scattered nodes layout where its numerical anisotropy~\cite{REUTHER201216} additionally distorts the model symmetry. Until now we assumed perfect symmetry and considered only envelope nodes from the first quadrant. Nevertheless, since a full dendrite is simulated we can also assess the level of anisotropy induced by the MIT scheme. For this purpose, the envelope nodes from all four quadrants are compared. On the left plot of Figure~\ref{fig:induced_anisotropy} the comparison of all four quadrants is presented for internodal spacing on the envelope set to $h_e \approx 0.028$ and at simulation time $t\approx 9.4$, i.e., close to the end of simulation, where the branches are well developed. All four envelope parts appear to be in good agreement.

For a more objective analysis, a cubic spline is fitted to the mapped nodes of the four envelope parts and the distance to the origin, $r_n(t) = \sqrt{x_n(t)^2 + y_n(t)^2}$, is measured with respect to the spline parameter, $t$, in the top tight corner in Figure~\ref{fig:induced_anisotropy}. Here, parameter $t$ runs from the top of the dendrite tip that grows along direction $\widehat{\b e}_y$ (vertically) towards the tip growing along direction $\widehat{\b e}_x$ (horizontally) and the index $n\in \{I, II, III, IV\}$ denotes the quadrants. The spline presentation allows us to evaluate the differences between the four curves, i.e., the four envelope parts from different quadrants.

The differences with respect to the envelope shape from the first quadrant, represented as a set of 1000 $(x_I, y_I)$ coordinates, are shown in Figure~\ref{fig:induced_anisotropy} in the bottom right. In this particular example, the largest difference is observed in the second quadrant, with the largest relative error slightly above $ 1\%$. The relative error in other quadrants is similar. The largest error, $e_\mathrm{max} = \max|(r_I-r_n)/r_I|$, is approximately 3 times the internodal distance enforced at the envelope, i.e., $e_\mathrm{max}/h_e \approx 3$. This may appear a lot, but note that it is only around $0.8\%$ of the grain size (given by the tip position $x_\mathrm{tip} \approx 8.88$).

%STo sum up, the anisotropy induced by the node arrangement is present, however, based on this analysis, we conclude that the contribution can be neglected.

\subsection{Spatial-discretization convergence}
Let us now take a closer look at the dendritic envelopes. To assess the behaviour of the numerical error in spatial discretization, we compare the envelopes computed with three different spatial discretizations at time $t\approx 9.6$, which corresponds to nearly the end of simulation time.

\begin{figure}[h]
	\centering
	\centering
	\includegraphics[width=0.8\textwidth]{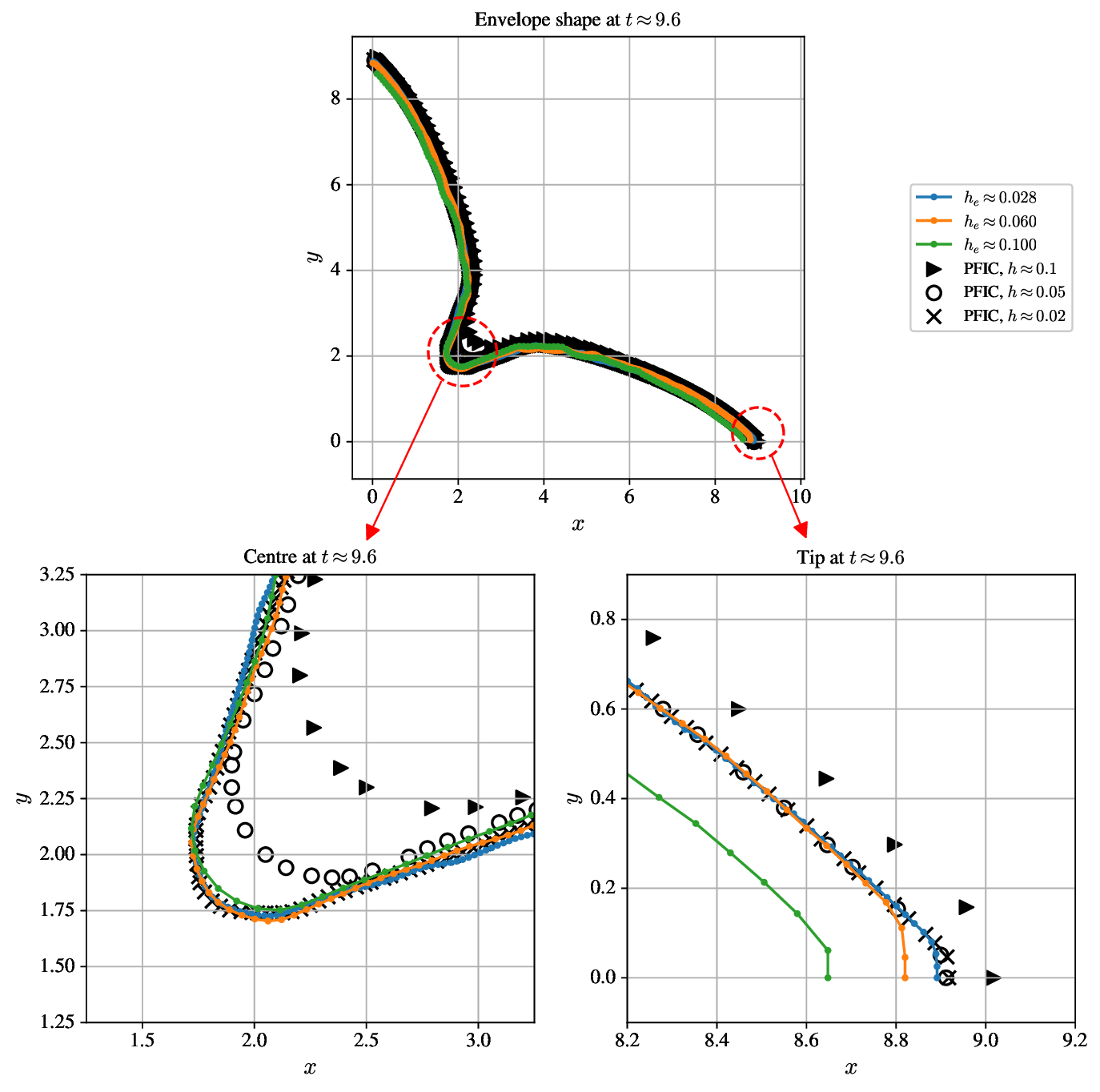}
	\caption{Convergence of the envelope shape at a selected simulation time $t\approx 9.6$. In the bottom row, the centre of the dendrite is shown on the left while the tip shape is shown on the right. \hl{For clarity, only every 20th node is plotted for the PFIC method.}}
	\label{fig:conv_shape}
\end{figure}

Figure~\ref{fig:conv_shape} shows the results for the finest ($\he \approx 0.028$), intermediate ($\he \approx 0.060$), and coarsest ($\he \approx 0.1$) using the proposed MIT-based approach and the results for the finest ($h \approx 0.02$), intermediate ($h \approx 0.05$), and coarsest ($h \approx 0.1$) discretizations using the PFIC-based approach. Note that the mesh spacing, $h$, of the PFIC is uniform and equal to $\he$ and that for the MIT, the pscing on the outer boundary, $h_\ell = 0.1$, is constant and identical for all simulations, regardless of $\he$.

\begin{figure}[H]
	\centering
	\includegraphics[width=0.85\textwidth]{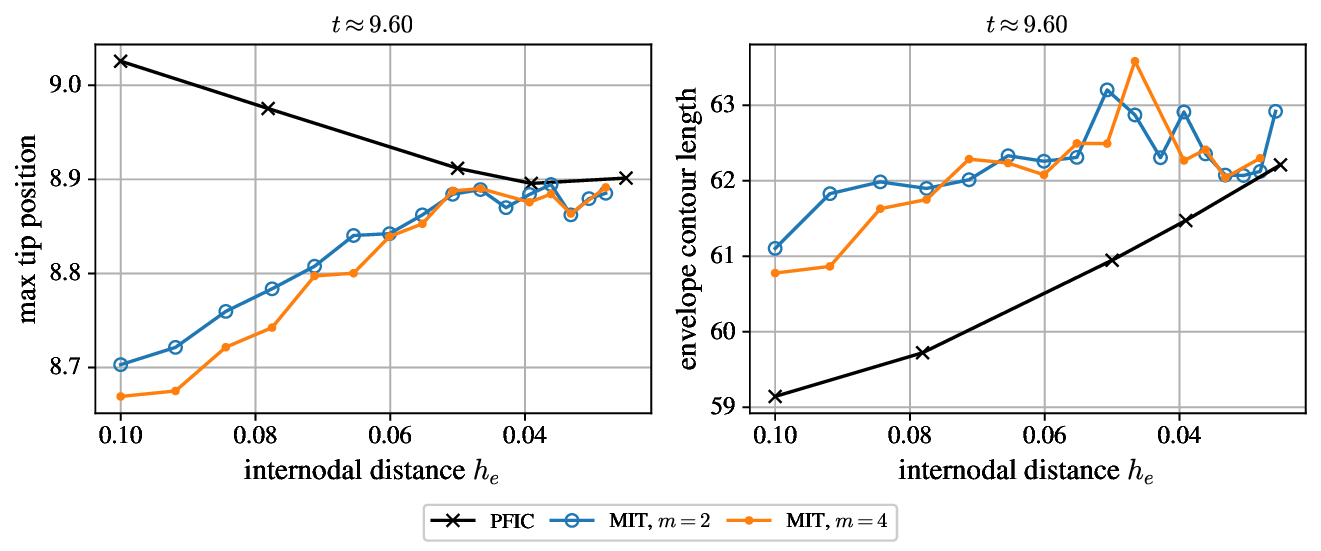}
	\caption{Maximum tip position with respect to the discretization spacing (left) and envelope contour length (right).}
	\label{fig:conv_envelope_length_pos}
\end{figure}

Convergent behaviour of the maximum tip position can be already seen in Figure~\ref{fig:conv_shape} at the bottom right figure.
%The maximum tip position computed with MIT is generally increasing with the number of the nodes on the envelope, while on the contrary the maximum tip positions computed with the PFIC-based approach decreases with respect to the quality of the domain discretization. 
To further support this observation, the maximum tip position is plotted with respect to the discretization spacing in Figure~\ref{fig:conv_envelope_length_pos} (left) near the end of the simulation. Here we can see that the slope of convergence of the MIT-based approach is steeper than that of PFIC and also that MIT reaches the asymptotic value at somewhat larger node spacing than PFIC. The PFIC-based approach converges towards slightly larger tip positions than those obtained with the MIT-based simulation. Nevertheless, the difference between the two is less than 0.5 \%, which is small, considering that we are comparing two conceptually different methods. We also show that increasing the order of the RBF approximation from second ($m=2$) to fourth ($m=4$) has only negligible impact on the solution of the full GEM. 

In the middle of the dendrite (see Figure~\ref{fig:conv_shape} in the lower left corner), we notice that MIT captures the final envelope shape with fewer nodes compared to PFIC. Note that PFIC with $\he \approx 0.1$ completely smooths out the depression in the centre, while MIT with a comparable spatial discretization captures all the details. The most likely reason for this is that MIT uses sharp interface tracking and avoids the smoothing effect of the phase field in curved parts of the shape, effectively requiring less nodes to reduce numerical artefacts. We attempt to evaluate this effect by measuring the length of the envelope. To do this, we use the spline representation of the dendritic envelope and calculate its total length by summing 5000 linear segments, which is shown in Figure~\ref{fig:conv_envelope_length_pos} (right).

\subsection{Computational cost}
\begin{figure}[h]
	\centering
	\centering
	\begin{subfigure}[b]{0.47\textwidth}
		\centering
		\includegraphics[width=\textwidth]{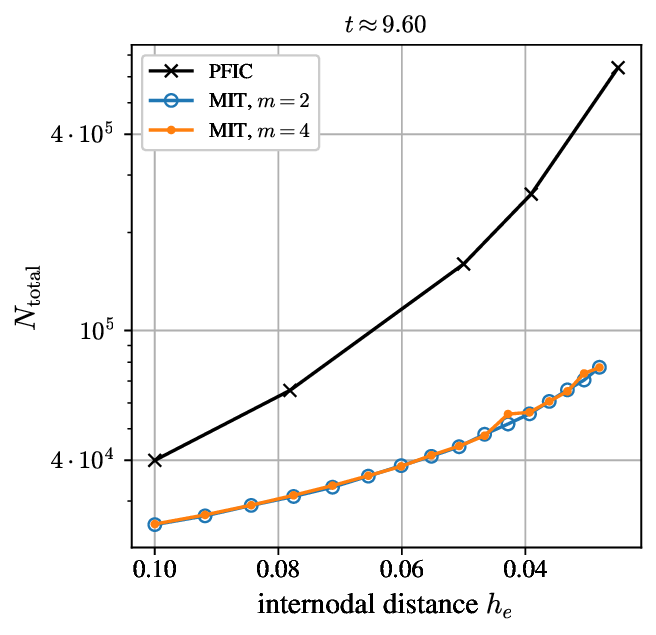}
		\caption{Number of nodes at simulation time $t\approx 9.60$ with respect to the internodal distance $\he$, comparison of MIT and PFIC.}
		\label{fig:node_count_at_time}
	\end{subfigure}
	\hfill
	\begin{subfigure}[b]{0.47\textwidth}
		\centering
		\includegraphics[width=\textwidth]{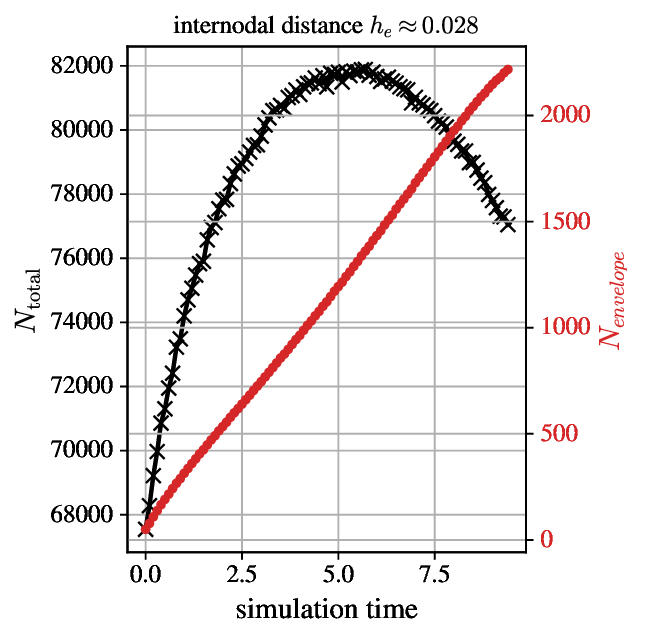}
		\caption{Time evolution of the number of nodes for the MIT for $\he=0.028$.}
		\label{fig:node_count_at_h}
	\end{subfigure}
	\caption{Analysis of the number of discretization nodes for \h-adaptive MIT and comparison to PFIC.}
	\label{fig:dendritic_envelopes2}
\end{figure}

To show the benefit of \h-adaptive discretization for the GEM we analyse the evolution of the number of discretization nodes during the simulations. Figure~\ref{fig:node_count_at_time} shows the total number of nodes (domain and boundary nodes), $N_\mathrm{total}$, with respect to the internodal distance $\he$. The number of nodes is generally much lower with the MIT approach, because in the MIT (i) the interior of the dendrite is not considered as part of the computational domain, and (ii) the $h$-refinement towards the dendrite envelope effectively leads to a lower total number of nodes, while preserving the quality of the local field description near the envelope. The last two remarks are particularly important for fine internodal distances $\he$, which yield about an order of magnitude fewer computational points for the MIT at $\he \approx 0.028$. Figure~\ref{fig:node_count_at_h} shows the evolution of the total number of discretization nodes with respect to the simulation time for the finest discretization used in this work ($\he \approx 0.028$). Given the \h-refinement, we find that the number of discretization points generally increases during the simulation. This is somewhat intuitive since the dendrite area increases with simulation time, i.e., the densely populated area increases. On the other hand, an increasing size of the dendrite also means a decrease in the total domain size, which prevails as the dendrite envelope approaches the outer boundary of the domain.

%\mz{Using space adaptivity (adaptive meshing or \h-adaptivity) in connection with the PFIC-GEM would allow a comparison of the impact of RBF-solution of diffusion and of the MIT inteface tracking on the performance.
%We did not implement space adaptivity to PFIC -- too complicated.
%We compare non-adaptive MIT to PFIC. This allows a very similar comparison because we omit the $h$-adaptivity from the both methods for the comparison. The comparison thus boils down to the benefit of RBF and of the MIT-interface tracking.
%
%Point 1: compares non-adaptive MIT to PFIC. The cost of the MIT is higher. We can suppose the the cost is:\\
%-decreased because of a smaller number of nodes (only the liquid part is in the domain for the solution of the diffusion equation)\\
%-decreased because it does not need to solve the PF equation for front propagation, but a computationally cheaper motion of markers.\\
%-decreased because of the different front tracking method that requires less evaluations of the front velocity ($O(N)$ instead of $O(9 N)$)\\
%-increased because of the node regeneration in the domain and on the envelope, required by the changing domain shape\\
%-increased because of the costlier RBF solution method for the diffusion equation\\
%The latter factors prevail and the cost of the MIT is around 50\% higher for a uniform node spacing.
%}

\begin{figure}
	\centering
	\includegraphics[width=0.7\textwidth]{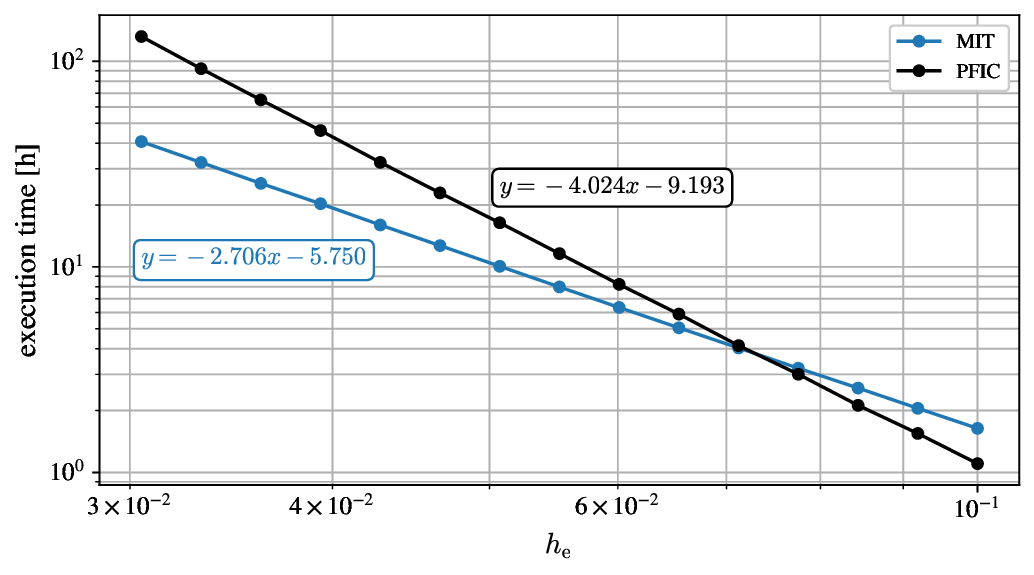}
	\caption{Execution times for PFIC and MIT methods as function of the node spacing at the envelope, $\he$. Note that for the MIT the node spacing at the outer boundary was constant ($h_\ell=0.1$) and for the PFIC the meshes are uniform, with spacing $\he$. The times were measured on the same hardware setup assuring CPU affinity and fixed CPU frequency at 2.4GHz. The regression functions give an estimation of the computational complexity.}
	\label{fig:execution_times}
\end{figure}

\hl{
To assess the computational cost of the MIT-GEM and to compare its performance to the PFIC-GEM we measured execution times in controlled conditions, i.e., on the same hardware setup, ensuring CPU affinity and fixed CPU frequency (2.4GHz). We ran MIT-GEM computations with $h_\ell=0.1$ on the outer boundary and with varying $\he$. The mesh spacing for the PFIC-GEM is uniform and equal to $\he$. Note that we can reasonably consider that the two methods have the same accuracy at equal $\he$. This is a simplification that is not far from reality (see Figure
}
~\ref{fig:conv_envelope_length_pos}). 
\hl{
The dependence of the execution time on $\he$ is shown in Figure
}
~\ref{fig:execution_times}.

\hl{
At $\he=h_\ell=0.1$ the node distribution of the MIT is uniform and this point therefore allows us to compare the performance for both methods at equivalent spatial discretization. In these conditions the execution time of the MIT-GEM is around 50\% longer. This is inherently due to the use of the RBF-FD method for the solution of the diffusion equation and due to the node regeneration required by the changing domain shape. These factors appear to prevail over those that are expected to decrease the computation time: smaller number of nodes, computationally cheaper front tracking by moving nodes instead of solving a phase-field equation, a smaller number of envelope velocity calculations.
}

\hl{
With refinement of the envelope node spacing, $\he$, we can observe a clear gain of the MIT-GEM over the PFIC-GEM. At the finest discretization, $\he=0.03$, the MIT-GEM is around three times faster. This gain is exclusively due to the \h-adaptivity of the spatial discretization in the MIT-GEM, which allows us to use significantly less nodes in the domain, compared to a uniform node distribution, while retaining a high node density on the envelope and thus a high accuracy of the solution in the regions of interest. For the problem at hand (the simulation of a single dendrite in a fairly large domain, constant node spacing at the outer boundary) the computational complexity of the MIT-GEM is shown in Figure
}
~\ref{fig:execution_times} 
\hl{
to scale as $\sim \he^{-2.7}$, where $\he$ is the node spacing on the envelope. For the PFIC-GEM it scales as $\sim \he^{-4}$, which is equivalent to $O(N^2)$, where $N$ is the number of mesh points. Note that the time step is $\sim \he^2$ for both methods. These results show that the MIT-GEM with \h-adaptivity is particularly interesting for high-accuracy solutions, using small $\he$, where \h-adaptivity is highly beneficial.
}

\section{Simulation of interacting grains}

\begin{figure}[H]
	\centering
	\includegraphics[width=\textwidth]{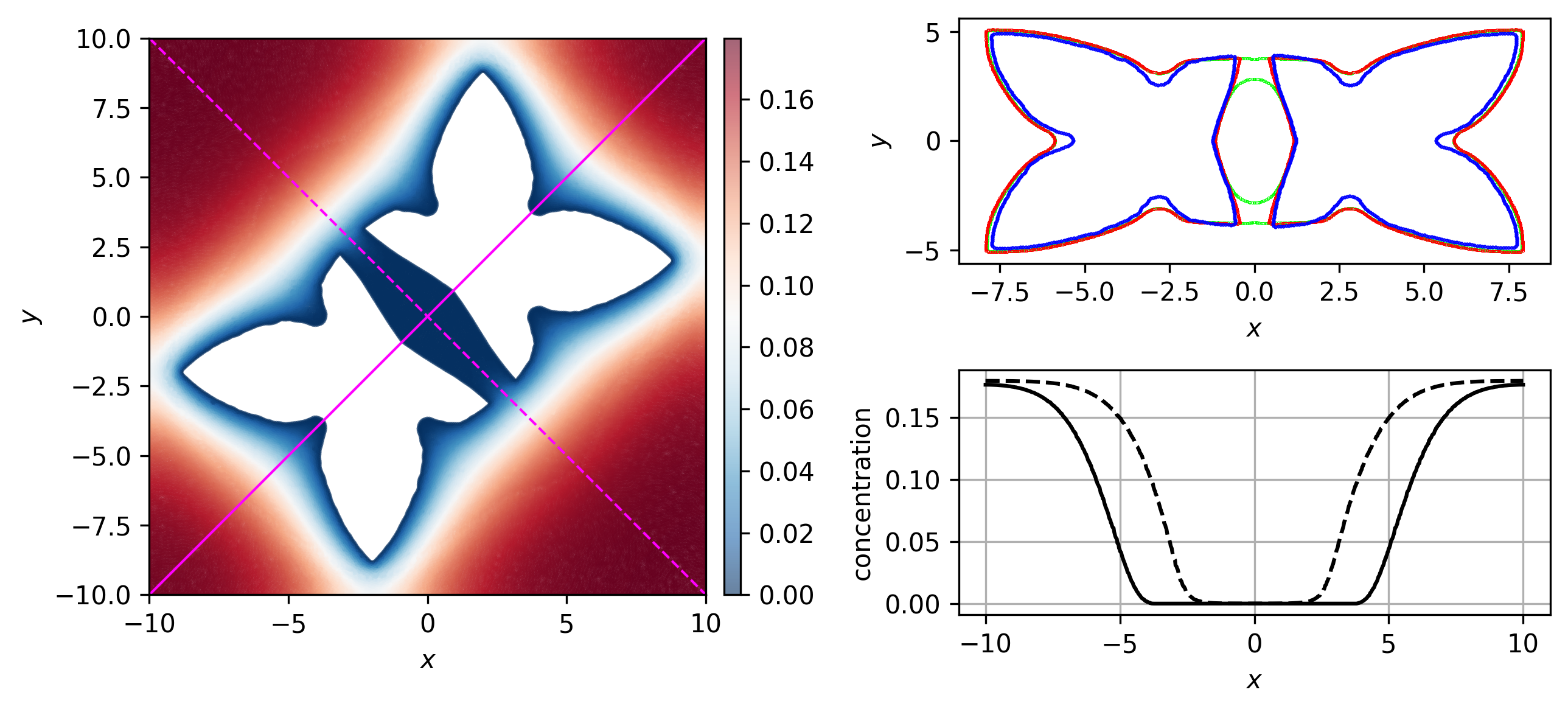}
	\caption{Two interacting dendrites at time $t=7.184$. Left: The MIT-GEM concentration field is shown. The envelopes predicted by PFIC-GEM using a single phase field are shown in top right by a green line, for two distinct phase fields with a red line while MIT envelopes are shown in blue. The grain seeds were positioned at $\b r_1 = (2, 2)$ and $\b r_2 = (-2, -2)$. The bottom right figure shows concentration profiles along the diagonals (shown in magenta in the left image).}
	\label{fig:concentration_double}
\end{figure}

\hl{
	One of the primary interests of the GEM is to study solutal interactions between grains over distances that are of the order of the grain size. In this section we show two examples of application of the MIT-GEM.

	In Figure
}
~\ref{fig:concentration_double}
\hl{
	we show MIT-GEM and PFIC-GEM simulations of two dendrites as they grow close to each other. A drawback of the PFIC when describing grains that are relatively close ($\sim \delta$) is that the phase-fields of the individual grains can interact and can thus provide envelope shapes that are not consistent with the GEM model. Due to the properties of the phase-field equations, two envelopes that are described by the same phase field can merge if they are sufficiently close. This is shown in Figure
}
~\ref{fig:concentration_double},
\hl{
	where two dendrite envelopes simulated by the PFIC-GEM using a single phase field  actually merge. This result is not consistent with the GEM model and the merging is only a consequence of the PFIC method. A way of avoiding merging of envelopes in the PFIC-GEM is to use a multi-phase-field approach, i.e., to use a distinct phase field for each grain. This is a robust method, but it adds an additional equation for each phase field, which increases the computation time and becomes costly when simulating large numbers of grains (100 or more). The MIT-GEM entirely avoids such merging due to numerical artefacts, the interactions between the envelopes are strictly as described by the GEM, they are given by the concentration field between the grains. The growth of the two envelopes shown in Figure
}
~\ref{fig:concentration_double}
\hl{
	stops at a distance of around 1.1. Note that the distance is somewhat larger than $\delta=1$ because the concentration of the liquid between the grains reaches 0 before the grains actually grow down to the distance of $\delta$. For clarity, we plot the solute concentration along the $y=x$ and $y=-x$ lines in Figure
}
~\ref{fig:concentration_double}
\hl{
	on the right. We see that the solute concentration between the dendrites is $0$, which stops further growth of the envelope in the MIT approach, consistently with the GEM model.
}

% \mz{Figure 18: Popravil sem caption. Dodati je treba še envelope za PFIC-GEM z dvema PF poljema (predlagam z rumeno črto, kot nakazano v captionu).}
%The simulation was started with two grains positioned at $\b r_1 = (2, 2)$ and $\b r_2 = (-2, -2)$ where the domain was discretized using $h_e\approx0.0775$ and $h_l=0.1$, and the time step was set to $\mathrm{d}t = 3.005\cdot 10^{-4}$.

\begin{figure}
	\centering
	\includegraphics[width=\textwidth]{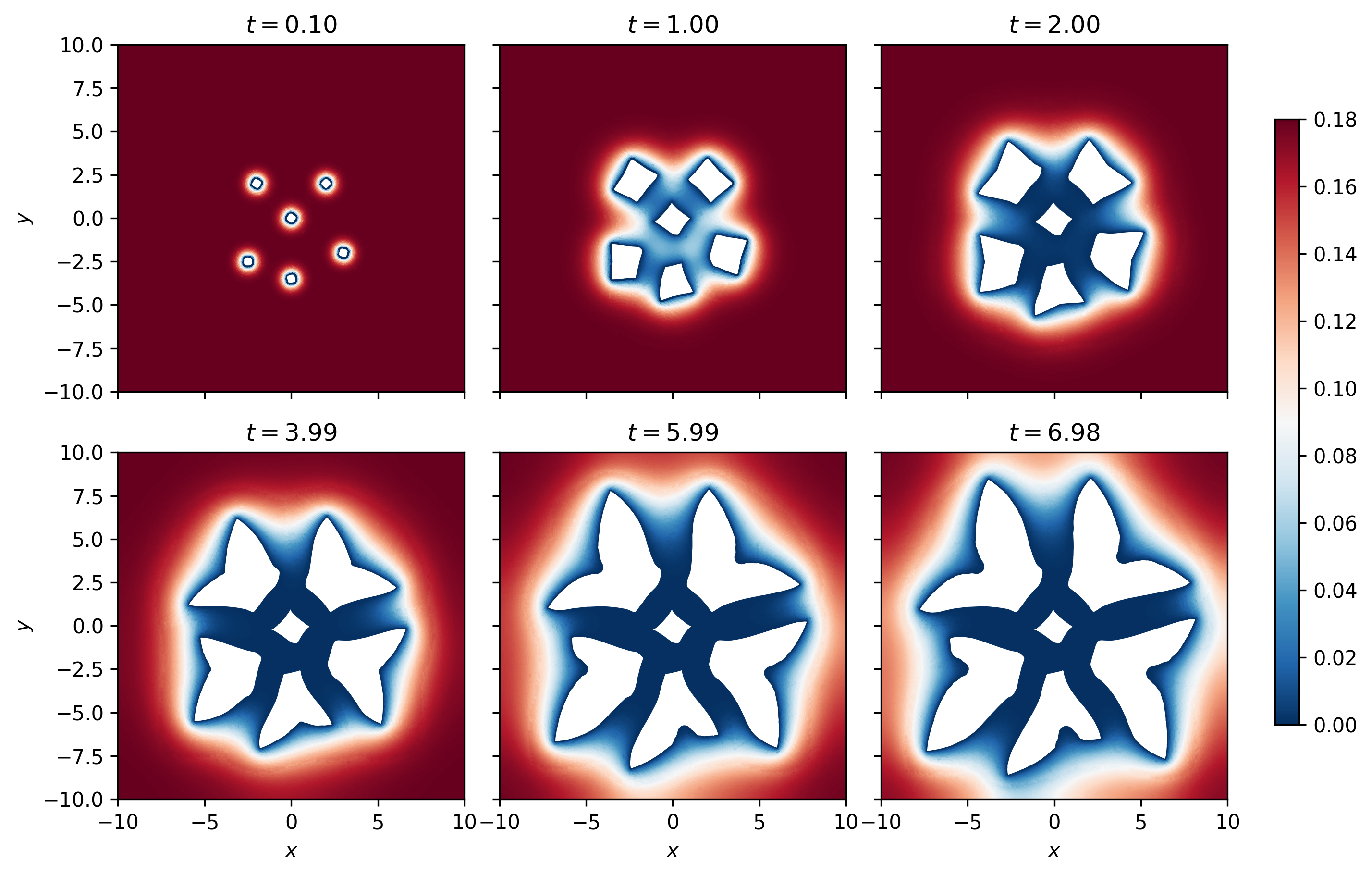}
	\caption{Growth of six interacting dendrites with different orientations of the tip growth directions.}
	\label{fig:time_evolution_6}
\end{figure}

\hl{
	To demonstrate a more general application of the MIT-GEM, we show a simulation of $6$ dendrites with different rotations of the four tip growth directions. The time evolution of the interacting grain shapes and of the solute concentration field is shown in Figure
}
~\ref{fig:time_evolution_6}.
\hl{
	Adding grains with different orientations in the MIT approach is straightforward. Each grain has its distinct set of envelope boundary nodes (evolving during the simulation) and can thus be easily attributed its own set of tip growth directions. This is an advantage with respect to the PFIC-GEM approach, where the multi-phase-field approach must be used to describe multiple grain orientations, introducing additional computational cost.
}

%
% Conclusions.
\section{Conclusions and perspectives}
\label{sec:conclusions}

In this paper we proposed the Meshless Interface Tracking (MIT) method for the solution of the Grain Envelope Model for dendritic growth (GEM). This novel method is expected to provide improved accuracy of GEM solutions in comparison with the Phase-Field Interface Capturing (PFIC) method, which has been used in all prior work on the GEM. Inherently, the MIT-GEM eliminates certain numerical artefacts experienced with the PFIC-GEM:
\begin{itemize}
	\item The MIT-GEM does not smooth the envelope, while smoothing cannot be entirely avoided with the PFIC-GEM.
	\item The MIT-GEM ensures that the envelope propagates exactly at the calculated speed, while the PFIC-GEM requires a careful surface reconstruction to ensure good accuracy of the envelope propagation.
\end{itemize}
Furthermore, the MIT introduces \h-adaptive spatial discretization, which is an advantage and enables one to ensure high accuracy of the solution of the diffusion field in the liquid. In the future it could also be used to refine the point spacing on the envelope in areas of high curvature. In addition to $h$ also an $hp$ adaptivity could be used~\cite{jancic_strong_2023}. 
\hl{
We analysed all critical elements of the MIT method, namely the meshless approximation of the differential operators, the interpolation of the concentration field, and the parametric envelope reconstruction used for rediscretization and for calculation of the envelope speed. We addressed the accuracy and computational complexity of the methods as well as the sensitivity to the order of RBF operator approximation and to the order of spline interpolation for envelope reconstruction. After the analysis we concluded that the order $k > 1$ should be used for surface reconstruction, and the second order RBF-FD approximation should be used for field reconstruction and its derivation. Note that these are not tuning parameters, but parameters that define the method (similar to the FEM, where the user decides whether to use linear basis functions or quadratic basis functions). Therefore, these parameters can be considered universal.
}

The comparison of PFIC and MIT shows that MIT can capture the shape of the envelope more accurately by using fewer nodes on the envelope. Moreover, due to the \h-refinement, the total number of nodes required in MIT is also much smaller compared to PFIC. MIT reaches high accuracy using second order methods for the approximation of the spatial partial differential operators, spatial field interpolation, and for spline reconstruction of the envelope. Higher order methods can be readily used for additional improvement of accuracy.

The focus of future applications of the MIT-GEM will be simulations requiring high accuracy, especially high accuracy of the envelope tracking. This includes simulations where the development of new branches of a dendrite is a key aspect. The MIT-GEM will also serve as a high-accuracy reference for the verification of the 
%computationally cheaper 
PFIC-GEM. From the numerical point of view, the introduction of a more complex node distribution that adapts to the local envelope curvature radius would be beneficial to improve the description of small-scale features of the envelope shape.

\hl{
The presented MIT is based on a meshless approximation using scattered nodes for spatial discretisation. Despite the 
%romantic 
often idealized
descriptions of meshless methods found in literature, there are still many unanswered questions that burden the method. The ongoing scientific debate about what is the best nodal distribution is inherently connected to the subsequent challenge of stencil selection, where scientific consensus is even more elusive. In MIT we use a widely accepted stencil that comprises the certain number of nearest nodes. Although such an approach is simple to implement and generalises to higher dimensions, it is also computationally expensive as relatively big stencils are required for stable approximation}~\cite{bayona2019insight, bayona2010rbf, bayona2017role}. 
\hl{The interplay of node positioning and stencil selection is of crucial importance not only for the execution performance but also for accuracy and stability, which is still, after decades of research, not fully understood}~\cite{kolar2023}. \hl{In the future work these concerns should be also tackled. In addition to the simple closest nodes stencils, there is promising research}~\cite{davydov2011adaptive,davydov2022improved}\hl{ showing that more sophisticated symmetric stencils could improve the stability and execution performance of the meshless methods. 

The introduction of scattered nodes also results in an unstructured nature of the data, which makes load balancing in parallel} execution~\cite{trobec_parallel_2015} and effective cache utilisation much more difficult~\cite{kosec_super_2014}, \hl{opening additional questions regarding the effective implementation of meshless methods, which could also be addressed in future work.}

% \newpage

% Acknowledgments
\section*{Acknowledgments}
M.J.\ and G.K.\ acknowledge the financial support from the Slovenian Research Agency research core funding No.\ P2-0095, and research projects No.\ J2-3048 and No. \ N2 - 0275. M.Z.\ acknowledges the support by the French State through the program ``Investment in the future'' operated by the National Research Agency (ANR) and referenced by ANR-11 LABX-0008-01 (LabEx DAMAS). A part of the required high performance computing resources was provided by the EXPLOR center hosted by the Université de Lorraine.

Funded by National Science Centre, Poland under the OPUS call in the Weave programme 2021/43/I/ST3/00228.
This research was funded in whole or in part by National Science Centre (2021/43/I/ST3/00228). For the purpose of Open Access,
the author has applied a CC-BY public copyright licence to any Author Accepted Manuscript (AAM) version arising from this submission.

\section*{Declarations}
\textbf{Conflict of interest.} The authors declare that they have no conflict of interest. All the co-authors have confirmed to know the submission of the manuscript by the corresponding author.

% References
\bibliographystyle{elsarticle-num}
\bibliography{references}

\begin{thebibliography}{10}
\expandafter\ifx\csname url\endcsname\relax
  \def\url#1{\texttt{#1}}\fi
\expandafter\ifx\csname urlprefix\endcsname\relax\def\urlprefix{URL }\fi
\expandafter\ifx\csname href\endcsname\relax
  \def\href#1#2{#2} \def\path#1{#1}\fi

\bibitem{steinbach1999three}
I.~Steinbach, C.~Beckermann, B.~Kauerauf, Q.~Li, J.~Guo, {Three-dimensional
  modeling of equiaxed dendritic growth on a mesoscopic scale}, Acta Materialia
  47 (1999) 971--982.
\newblock \href {https://doi.org/10.1016/S1359-6454(98)00380-2}
  {\path{doi:10.1016/S1359-6454(98)00380-2}}.

\bibitem{delaleau2010mesoscopic}
P.~Delaleau, C.~Beckermann, R.~H. Mathiesen, L.~Arnberg, {Mesoscopic Simulation
  of Dendritic Growth Observed in X-ray Video Microscopy During Directional
  Solidification of Al–Cu Alloys}, ISIJ International 50 (2010) 1886--1894.
\newblock \href {https://doi.org/10.2355/isijinternational.50.1886}
  {\path{doi:10.2355/isijinternational.50.1886}}.

\bibitem{souhar2016three}
Y.~Souhar, V.~F. De~Felice, C.~Beckermann, H.~Combeau, M.~Zalo{\v{z}}nik,
  Three-dimensional mesoscopic modeling of equiaxed dendritic solidification of
  a binary alloy, Computational Materials Science 112 (2016) 304--317.

\bibitem{steinbach2005transient}
I.~Steinbach, H.-J. Diepers, C.~Beckermann, Transient growth and interaction of
  equiaxed dendrites, Journal of crystal growth 275~(3-4) (2005) 624--638.

\bibitem{olmedilla2019quantitative}
A.~Olmedilla, M.~Zalo{\v{z}}nik, H.~Combeau, {Quantitative 3D mesoscopic
  modeling of grain interactions during equiaxed dendritic solidification in a
  thin sample}, Acta Materialia 173 (2019) 249--261.
\newblock \href {https://doi.org/10.1016/j.actamat.2019.05.019}
  {\path{doi:10.1016/j.actamat.2019.05.019}}.

\bibitem{chirouf2023investigation}
A.~Chirouf, B.~Appolaire, A.~Finel, Y.~{Le Bouar}, M.~Zalo{\v{z}}nik,
  \href{https://dx.doi.org/10.1088/1757-899X/1281/1/012054}{Investigation of
  diffusive grain interactions during equiaxed dendritic solidification}, IOP
  Conference Series: Materials Science and Engineering 1281 (2023) 012054.
\newblock \href {https://doi.org/10.1088/1757-899X/1281/1/012054}
  {\path{doi:10.1088/1757-899X/1281/1/012054}}.
\newline\urlprefix\url{https://dx.doi.org/10.1088/1757-899X/1281/1/012054}

\bibitem{viardin2017mesoscopic}
A.~Viardin, M.~Zalo{\v{z}}nik, Y.~Souhar, M.~Apel, H.~Combeau, {Mesoscopic
  modeling of spacing and grain selection in columnar dendritic solidification:
  Envelope versus phase-field model}, Acta Materialia 122 (2017) 386--399.
\newblock \href {https://doi.org/10.1016/j.actamat.2016.10.004}
  {\path{doi:10.1016/j.actamat.2016.10.004}}.

\bibitem{viardin2020mesoscopic}
A.~Viardin, Y.~Souhar, M.~{Cisternas Fern{\'{a}}ndez}, M.~Apel,
  M.~Zalo{\v{z}}nik, {Mesoscopic modeling of equiaxed and columnar
  solidification microstructures under forced flow and buoyancy-driven flow in
  hypergravity: Envelope versus phase-field model}, Acta Materialia 199 (2020)
  680--694.
\newblock \href {https://doi.org/10.1016/j.actamat.2020.07.069}
  {\path{doi:10.1016/j.actamat.2020.07.069}}.

\bibitem{torabirad2019upscaling}
M.~{Torabi Rad}, M.~Zalo{\v{z}}nik, H.~Combeau, C.~Beckermann, {Upscaling
  mesoscopic simulation results to develop constitutive relations for
  macroscopic modeling of equiaxed dendritic solidification}, Materialia
  5~(November 2018) (2019).
\newblock \href {https://doi.org/10.1016/j.mtla.2019.100231}
  {\path{doi:10.1016/j.mtla.2019.100231}}.

\bibitem{tourret2020comparing}
D.~Tourret, L.~Sturz, A.~Viardin, M.~Zalo{\v{z}}nik, Comparing mesoscopic
  models for dendritic growth, in: IOP Conference Series: Materials Science and
  Engineering, Vol. 861, IOP Publishing, 2020, p. 012002.

\bibitem{boukellal2023multiscale}
A.~K. Boukellal, M.~Zalo{\v{z}}nik, J.-M. Debierre,
  \href{https://dx.doi.org/10.1088/1757-899X/1281/1/012048}{Multi-scale
  modeling of equiaxed dendritic solidification of al-cu at constant cooling
  rate}, IOP Conference Series: Materials Science and Engineering 1281 (2023)
  012048.
\newblock \href {https://doi.org/10.1088/1757-899X/1281/1/012048}
  {\path{doi:10.1088/1757-899X/1281/1/012048}}.
\newline\urlprefix\url{https://dx.doi.org/10.1088/1757-899X/1281/1/012048}

\bibitem{sun2007sharp}
Y.~Sun, C.~Beckermann, sharp-interface tracking using the phase-field equation,
  Journal of Computational Physics 220~(2) (2007) 626--653.

\bibitem{Zhang2019d}
A.~Zhang, J.~Du, Z.~Guo, Q.~Wang, S.~Xiong, {Conservative phase-field method
  with a parallel and adaptive-mesh-refinement technique for interface
  tracking}, Physical Review E 100~(2) (2019) 23305.
\newblock \href {https://doi.org/10.1103/PhysRevE.100.023305}
  {\path{doi:10.1103/PhysRevE.100.023305}}.

\bibitem{Provatas1999}
N.~Provatas, N.~Goldenfeld, J.~Dantzig, {Adaptive Mesh Refinement Computation
  of Solidification Microstructures Using Dynamic Data Structures}, Journal of
  Computational Physics 148~(1) (1999) 265--290.
\newblock \href {http://arxiv.org/abs/9808216} {\path{arXiv:9808216}}, \href
  {https://doi.org/10.1006/jcph.1998.6122} {\path{doi:10.1006/jcph.1998.6122}}.

\bibitem{Sarkis2016}
C.~Sarkis, L.~Silva, C.~A. Gandin, M.~Plapp, {Three-dimensional modeling of a
  thermal dendrite using the phase field method with automatic anisotropic and
  unstructured adaptive finite element meshing}, IOP Conference Series:
  Materials Science and Engineering 117~(1) (2016).
\newblock \href {https://doi.org/10.1088/1757-899X/117/1/012008}
  {\path{doi:10.1088/1757-899X/117/1/012008}}.

\bibitem{Greenwood2018}
M.~Greenwood, K.~N. Shampur, N.~Ofori-Opoku, T.~Pinomaa, L.~Wang, S.~Gurevich,
  N.~Provatas, {Quantitative 3D phase field modelling of solidification using
  next-generation adaptive mesh refinement}, Computational Materials Science
  142 (2018) 153--171.
\newblock \href {https://doi.org/10.1016/j.commatsci.2017.09.029}
  {\path{doi:10.1016/j.commatsci.2017.09.029}}.

\bibitem{Sakane2022}
S.~Sakane, T.~Aoki, T.~Takaki, {Parallel GPU-accelerated adaptive mesh
  refinement on two-dimensional phase-field lattice Boltzmann simulation of
  dendrite growth}, Computational Materials Science 211~(May) (2022) 111507.
\newblock \href {https://doi.org/10.1016/j.commatsci.2022.111507}
  {\path{doi:10.1016/j.commatsci.2022.111507}}.

\bibitem{Guo2015a}
Z.~Guo, S.~M. Xiong, {On solving the 3-D phase field equations by employing a
  parallel-adaptive mesh refinement (Para-AMR) algorithm}, Computer Physics
  Communications 190 (2015) 89--97.
\newblock \href {https://doi.org/10.1016/j.cpc.2015.01.016}
  {\path{doi:10.1016/j.cpc.2015.01.016}}.

\bibitem{Ham2024}
S.~Ham, Y.~Li, S.~Kwak, D.~Jeong, J.~Kim, {An efficient and fast adaptive
  numerical method for a novel phase-field model of crystal growth},
  Communications in Nonlinear Science and Numerical Simulation 131~(November
  2023) (2024) 107822.
\newblock \href {https://doi.org/10.1016/j.cnsns.2024.107822}
  {\path{doi:10.1016/j.cnsns.2024.107822}}.

\bibitem{Dobravec2022}
T.~Dobravec, B.~Mavri{\v{c}}, B.~{\v{S}}arler, {Acceleration of RBF-FD meshless
  phase-field modelling of dendritic solidification by space-time adaptive
  approach}, Computers and Mathematics with Applications 126~(September) (2022)
  77--99.
\newblock \href {https://doi.org/10.1016/j.camwa.2022.09.008}
  {\path{doi:10.1016/j.camwa.2022.09.008}}.

\bibitem{ghoneim2020smoothed}
A.~Y. Ghoneim, A smoothed particle hydrodynamics-phase field method with radial
  basis functions and moving least squares for meshfree simulation of dendritic
  solidification, Applied Mathematical Modelling 77 (2020) 1704--1741.

\bibitem{bahramifar2022local}
S.~Bahramifar, F.~Mossaiby, H.~Haftbaradaran, A local meshless method for
  transient nonlinear problems: Preliminary investigation and application to
  phase-field models, Computers \& Mathematics with Applications 124 (2022)
  163--187.

\bibitem{ghoneim2016new}
A.~Ghoneim, A new technique for numerical simulation of dendritic
  solidification using a meshfree interface finite element method,
  International Journal for Numerical Methods in Engineering 107~(10) (2016)
  813--852.

\bibitem{Gibou2003}
F.~Gibou, R.~Fedkiw, R.~Caflisch, S.~Osher, {A Level Set Approach for the
  Numerical Simulation of Dendritic Growth}, Journal of Scientific Computing
  19~(1-3) (2003) 183--199.
\newblock \href {https://doi.org/10.1023/A:1025399807998}
  {\path{doi:10.1023/A:1025399807998}}.

\bibitem{Tan2006}
L.~Tan, N.~Zabaras, {A level set simulation of dendritic solidification with
  combined features of front-tracking and fixed-domain methods}, Journal of
  Computational Physics 211~(1) (2006) 36--63.
\newblock \href {https://doi.org/10.1016/j.jcp.2005.05.013}
  {\path{doi:10.1016/j.jcp.2005.05.013}}.

\bibitem{Ghoneim2018}
A.~Y. Ghoneim, {The Meshfree Interface Finite Element Method for Numerical
  Simulation of Dendritic Solidification with Fluid Flow}, International
  Journal of Computational Methods 15~(7) (2018).
\newblock \href {https://doi.org/10.1142/S0219876218500573}
  {\path{doi:10.1142/S0219876218500573}}.

\bibitem{Ramanuj2019}
V.~Ramanuj, R.~Sankaran, B.~Radhakrishnan, {A sharp interface model for
  deterministic simulation of dendrite growth}, Computational Materials Science
  169~(May) (2019) 109097.
\newblock \href {https://doi.org/10.1016/j.commatsci.2019.109097}
  {\path{doi:10.1016/j.commatsci.2019.109097}}.

\bibitem{Limare2023}
A.~Limare, S.~Popinet, C.~Josserand, Z.~Xue, A.~Ghigo, {A hybrid level-set /
  embedded boundary method applied to solidification-melt problems}, Journal of
  Computational Physics 474 (2023).
\newblock \href {https://doi.org/10.1016/j.jcp.2022.111829}
  {\path{doi:10.1016/j.jcp.2022.111829}}.

\bibitem{DuChene2008}
A.~{Du Ch{\'{e}}n{\'{e}}}, C.~Min, F.~Gibou, {Second-order accurate computation
  of curvatures in a level set framework using novel high-order
  reinitialization schemes}, Journal of Scientific Computing 35~(2-3) (2008)
  114--131.
\newblock \href {https://doi.org/10.1007/s10915-007-9177-1}
  {\path{doi:10.1007/s10915-007-9177-1}}.

\bibitem{Juric1996}
D.~Juric, G.~Tryggvason, {A front-tracking method for dendritic
  solidification}, Journal of Computational Physics 123~(1) (1996) 127--148.
\newblock \href {https://doi.org/10.1006/jcph.1996.0011}
  {\path{doi:10.1006/jcph.1996.0011}}.

\bibitem{Reuther2014a}
K.~Reuther, M.~Rettenmayr, {Simulating dendritic solidification using an
  anisotropy-free meshless front-tracking method}, Journal of Computational
  Physics 279 (2014) 63--66.
\newblock \href {https://doi.org/10.1016/j.jcp.2014.09.003}
  {\path{doi:10.1016/j.jcp.2014.09.003}}.

\bibitem{nguyen_meshless_2008}
V.~P. Nguyen, T.~Rabczuk, S.~Bordas, M.~Duflot,
  \href{https://sourceforge.net/projects/elemfregalerkin/}{Meshless methods: a
  review and computer implementation aspects}, Math. Comput. Simul 79~(3)
  (2008) 763--813.
\newblock \href {https://doi.org/10.1016/j.matcom.2008.01.003}
  {\path{doi:10.1016/j.matcom.2008.01.003}}.
\newline\urlprefix\url{https://sourceforge.net/projects/elemfregalerkin/}

\bibitem{bayona2019insight}
V.~Bayona, An insight into rbf-fd approximations augmented with polynomials,
  Computers \& Mathematics with Applications 77~(9) (2019) 2337--2353.

\bibitem{tolstykh2003using}
A.~Tolstykh, D.~Shirobokov, On using radial basis functions in a “finite
  difference mode” with applications to elasticity problems, Computational
  Mechanics 33~(1) (2003) 68--79.

\bibitem{slak2019generation}
J.~Slak, G.~Kosec, On generation of node distributions for meshless pde
  discretizations, SIAM journal on scientific computing 41~(5) (2019)
  A3202--A3229.

\bibitem{ortega_meshless_2013}
E.~Ortega, E.~Oñate, S.~Idelsohn, R.~Flores, A meshless finite point method
  for three-dimensional analysis of compressible flow problems involving moving
  boundaries and adaptivity, International Journal for Numerical Methods in
  Fluids 73~(4) (2013) 323--343, publisher: Wiley Online Library.

\bibitem{slak2019adaptive}
J.~Slak, G.~Kosec, Adaptive radial basis function--generated finite differences
  method for contact problems, International Journal for Numerical Methods in
  Engineering 119~(7) (2019) 661--686.

\bibitem{cantor1977dendritic}
B.~Cantor, A.~Vogel, Dendritic solidification and fluid flow, Journal of
  Crystal Growth 41~(1) (1977) 109--123.

\bibitem{sun2008atwo}
Y.~Sun, C.~Beckermann, {A two-phase diffuse-interface model for Hele–Shaw
  flows with large property contrasts}, Physica D: Nonlinear Phenomena 237
  (2008) 3089--3098.
\newblock \href {https://doi.org/10.1016/j.physd.2008.06.010}
  {\path{doi:10.1016/j.physd.2008.06.010}}.

\bibitem{le2023guidelines}
S.~Le~Borne, W.~Leinen, Guidelines for rbf-fd discretization: Numerical
  experiments on the interplay of a multitude of parameter choices, Journal of
  Scientific Computing 95~(1) (2023) 8.

\bibitem{tominec2021least}
I.~Tominec, E.~Larsson, A.~Heryudono, A least squares radial basis function
  finite difference method with improved stability properties, SIAM Journal on
  Scientific Computing 43~(2) (2021) A1441--A1471.

\bibitem{jancic_monomial_2021}
M.~Jančič, J.~Slak, G.~Kosec, Monomial {Augmentation} {Guidelines} for
  {RBF}-{FD} from {Accuracy} {Versus} {Computational} {Time} {Perspective},
  Journal of Scientific Computing 87~(1), publisher: Springer Science and
  Business Media LLC (Feb. 2021).
\newblock \href {https://doi.org/10.1007/s10915-020-01401-y}
  {\path{doi:10.1007/s10915-020-01401-y}}.

\bibitem{depolli_parallel_2022}
M.~Depolli, J.~Slak, G.~Kosec,
  \href{https://linkinghub.elsevier.com/retrieve/pii/S0045794922000335}{Parallel
  domain discretization algorithm for {RBF}-{FD} and other meshless numerical
  methods for solving {PDEs}}, Computers \& Structures 264 (2022) 106773.
\newblock \href {https://doi.org/10.1016/j.compstruc.2022.106773}
  {\path{doi:10.1016/j.compstruc.2022.106773}}.
\newline\urlprefix\url{https://linkinghub.elsevier.com/retrieve/pii/S0045794922000335}

\bibitem{de2019fast}
S.~De~Marchi, A.~Mart{\'\i}nez, E.~Perracchione, Fast and stable rational
  rbf-based partition of unity interpolation, Journal of Computational and
  Applied Mathematics 349 (2019) 331--343.

\bibitem{jancic-stability}
M.~Jančič, G.~Kosec, Stability analysis of rbf-fd and wls based local strong
  form meshless methods on scattered nodes, in: 2022 45th Jubilee International
  Convention on Information, Communication and Electronic Technology (MIPRO),
  2022, pp. 275--280.
\newblock \href {https://doi.org/10.23919/MIPRO55190.2022.9803334}
  {\path{doi:10.23919/MIPRO55190.2022.9803334}}.

\bibitem{davydov2022improved}
O.~Davydov, D.~T. Oanh, N.~M. Tuong, Improved stencil selection for meshless
  finite difference methods in 3d, Journal of Computational and Applied
  Mathematics (2023) 115031.

\bibitem{9803651}
M.~Rot, A.~Rashkovska, Meshless method stencil evaluation with machine
  learning, in: 2022 45th Jubilee International Convention on Information,
  Communication and Electronic Technology (MIPRO), 2022, pp. 269--274.
\newblock \href {https://doi.org/10.23919/MIPRO55190.2022.9803651}
  {\path{doi:10.23919/MIPRO55190.2022.9803651}}.

\bibitem{bayona2017role}
V.~Bayona, N.~Flyer, B.~Fornberg, G.~A. Barnett, On the role of polynomials in
  rbf-fd approximations: Ii. numerical solution of elliptic pdes, Journal of
  Computational Physics 332 (2017) 257--273.

\bibitem{flyer2016role}
N.~Flyer, B.~Fornberg, V.~Bayona, G.~A. Barnett, On the role of polynomials in
  rbf-fd approximations: I. interpolation and accuracy, Journal of
  Computational Physics 321 (2016) 21--38.

\bibitem{pu-interpolation}
J.~Slak, Partition-of-unity based error indicator for local collocation
  meshless methods, in: 2021 44th International Convention on Information,
  Communication and Electronic Technology (MIPRO), 2021, pp. 254--258.
\newblock \href {https://doi.org/10.23919/MIPRO52101.2021.9597066}
  {\path{doi:10.23919/MIPRO52101.2021.9597066}}.

\bibitem{DEMARCHI2019331}
S.~{De Marchi}, A.~Martínez, E.~Perracchione,
  \href{https://www.sciencedirect.com/science/article/pii/S0377042718304400}{Fast
  and stable rational rbf-based partition of unity interpolation}, Journal of
  Computational and Applied Mathematics 349 (2019) 331--343.
\newblock \href {https://doi.org/https://doi.org/10.1016/j.cam.2018.07.020}
  {\path{doi:https://doi.org/10.1016/j.cam.2018.07.020}}.
\newline\urlprefix\url{https://www.sciencedirect.com/science/article/pii/S0377042718304400}

\bibitem{wendland2004scattered}
H.~Wendland, Scattered data approximation, Vol.~17, Cambridge university press,
  2004.

\bibitem{bayona2010rbf}
V.~Bayona, M.~Moscoso, M.~Carretero, M.~Kindelan, Rbf-fd formulas and
  convergence properties, Journal of Computational Physics 229~(22) (2010)
  8281--8295.

\bibitem{surface}
M.~Jančič, V.~Cvrtila, G.~Kosec, Discretized boundary surface reconstruction,
  in: 2021 44th International Convention on Information, Communication and
  Electronic Technology (MIPRO), 2021, pp. 278--283.
\newblock \href {https://doi.org/10.23919/MIPRO52101.2021.9596965}
  {\path{doi:10.23919/MIPRO52101.2021.9596965}}.

\bibitem{duh_surface}
U.~Duh, G.~Kosec, J.~Slak, Fast variable density node generation on parametric
  surfaces with application to mesh-free methods, SIAM Journal on Scientific
  Computing 43~(2) (2021) A980--A1000.
\newblock \href {https://doi.org/10.1137/20M1325642}
  {\path{doi:10.1137/20M1325642}}.

\bibitem{slak2021medusa}
J.~Slak, G.~Kosec, Medusa: A c++ library for solving pdes using strong form
  mesh-free methods, ACM Transactions on Mathematical Software (TOMS) 47~(3)
  (2021) 1--25.

\bibitem{Dobravec2023}
T.~Dobravec, B.~Mavri{\v{c}}, R.~Zahoor, B.~{\v{S}}arler, {A coupled
  domain–boundary type meshless method for phase-field modelling of dendritic
  solidification with the fluid flow}, International Journal of Numerical
  Methods for Heat and Fluid Flow 33~(8) (2023) 2963--2981.
\newblock \href {https://doi.org/10.1108/HFF-03-2023-0131}
  {\path{doi:10.1108/HFF-03-2023-0131}}.

\bibitem{REUTHER201216}
K.~Reuther, B.~Sarler, M.~Rettenmayr,
  \href{https://www.sciencedirect.com/science/article/pii/S1290072911002559}{Solving
  diffusion problems on an unstructured, amorphous grid by a meshless method},
  International Journal of Thermal Sciences 51 (2012) 16--22.
\newblock \href
  {https://doi.org/https://doi.org/10.1016/j.ijthermalsci.2011.08.017}
  {\path{doi:https://doi.org/10.1016/j.ijthermalsci.2011.08.017}}.
\newline\urlprefix\url{https://www.sciencedirect.com/science/article/pii/S1290072911002559}

\bibitem{jancic_strong_2023}
M.~Jančič, G.~Kosec,
  \href{https://link.springer.com/10.1007/s00366-023-01843-6}{Strong form
  mesh-free hp-adaptive solution of linear elasticity problem}, Engineering
  with Computers (May 2023).
\newblock \href {https://doi.org/10.1007/s00366-023-01843-6}
  {\path{doi:10.1007/s00366-023-01843-6}}.
\newline\urlprefix\url{https://link.springer.com/10.1007/s00366-023-01843-6}

\bibitem{kolar2023}
A.~Kolar-Požun, M.~Jančič, M.~Rot, G.~Kosec, Oscillatory behaviour of the
  rbf-fd approximation accuracy under increasing stencil size, in:
  Computational Science - ICCS 2023 : 23rd International Conference, 2023.

\bibitem{davydov2011adaptive}
O.~Davydov, D.~T. Oanh, Adaptive meshless centres and rbf stencils for poisson
  equation, Journal of Computational Physics 230~(2) (2011) 287--304.

\bibitem{trobec_parallel_2015}
R.~Trobec, G.~Kosec, Parallel scientific computing: theory, algorithms, and
  applications of mesh based and meshless methods, Springer, 2015.

\bibitem{kosec_super_2014}
G.~Kosec, M.~Depolli, A.~Rashkovska, R.~Trobec, Super linear speedup in a local
  parallel meshless solution of thermo-fluid problems, Computers \& Structures
  133 (2014) 30--38.
\newblock \href {https://doi.org/10.1016/j.compstruc.2013.11.016}
  {\path{doi:10.1016/j.compstruc.2013.11.016}}.

\end{thebibliography}

\end{document}